\documentclass[11pt]{article}
\usepackage[short]{optidef}
\usepackage{tikz, lscape, amsmath,amsfonts,latexsym,amssymb, amsthm}
\usepackage[utf8]{inputenc}
\usepackage[english]{babel}
\usepackage{lmodern, verbatim}
\usepackage[T1]{fontenc}
\usepackage{calc,bm}
\usepackage{tabularx,enumitem,natbib,nameref,hyperref}
\usepackage[mathscr]{eucal}
\usepackage{subfig,stackrel,floatrow}

\renewcommand{\b}{\mathbf}

\usepackage[margin=1in]{geometry}

\newcommand{\BB}{\mathbb{B}}

\newcommand{\epsArr}{\overrightarrow{\bm{\epsilon}}}
\newcommand{\eArr}{\overrightarrow{\bm{e}}}
\newcommand{\etaArr}{\overrightarrow{\bm{\eta}}}
\newcommand{\omeArr}{\overrightarrow{\bm{\omega}}}
\newcommand{\NN}{\mathbb{N}}

\newcommand{\RR}{\mathbb{R}}

\renewcommand{\H}{\mathcal{H}}

\newcommand{\Hk}{\mathcal{H}_K}

\newcommand{\Psc}{\mathscr{P}}
\newcommand{\Lcal}{\mathcal{L}}
\newcommand{\Sx}{\mathcal{S}}

\newcommand{\A}{\mathcal{A}}
\newcommand{\V}{\mathcal{V}}
\newcommand{\C}{\mathcal{C}}

\newcommand{\T}{\mathscr{T}}

\DeclareMathOperator*{\Sp}{span}

\newcommand{\Id}{\text{Id}}
\renewcommand{\d}{\mathrm{d}} 
\newcommand{\tb}{\textbf} 
\newcommand{\R}{\RR} 
\renewcommand{\b}{\mathbf}
\newcommand{\rv}{}
\newcommand{\rvS}{}
 
\newcommand{\iv}[2]{[\![#1,#2]\!]} 

\theoremstyle{definition}

\theoremstyle{plain}
\newtheorem{Theorem}{Theorem}

\newtheorem{Proposition}{Proposition}
\newtheorem{Lemma}{Lemma}

\theoremstyle{definition}
\newtheorem{Definition}{Definition}
\date{\today}

\makeatletter
\let\orgdescriptionlabel\descriptionlabel
\renewcommand*{\descriptionlabel}[1]{%
	\let\orglabel\label
	\let\label\@gobble
	\phantomsection
	\edef\@currentlabel{#1\unskip}%
	\let\label\orglabel
	\orgdescriptionlabel{#1}%
}
\makeatother

\title{Linearly-constrained Linear Quadratic Regulator from the viewpoint of kernel methods\footnote{To be published in 2021 in SIAM Journal on Control and Optimization}}
\author{Pierre-Cyril Aubin-Frankowski\footnote{\'Ecole des Ponts ParisTech and CAS, MINES ParisTech, PSL Research University, France. Email: pierre-cyril.aubin@mines-paristech.fr}}
\date{\today}

\begin{document}
	\maketitle
	\begin{abstract}
		The linear quadratic regulator problem is central in optimal control and was investigated since the very beginning of control theory. Nevertheless, when it includes affine state constraints, it remains very challenging from the classical ``maximum principle`` perspective. In this study we present how matrix-valued reproducing kernels allow for an alternative viewpoint. We show that the quadratic objective paired with the linear dynamics encode the relevant kernel, defining a Hilbert space of controlled trajectories. Drawing upon kernel formalism, we introduce a strengthened continuous-time convex optimization problem which can be tackled exactly with finite dimensional solvers, and which solution is interior to the constraints. When refining a time-discretization grid, this solution can be made arbitrarily close to the solution of the state-constrained Linear Quadratic Regulator. We illustrate the implementation of this method on a path-planning problem.
	\end{abstract}

	\section{Introduction}\label{secIntro}
	
	In its simplest form, the problem of time-varying linear quadratic optimal control with finite horizon and affine inequality state constraints writes as
	\begin{mini}
		{\substack{\b x(\cdot),\b u(\cdot)}}{g(\b x(T))+\int_{0}^{T}[\b x(t)^{\top}\b Q(t) \b x(t) + \b u(t)^{\top} \b R(t) \b u(t)] \d t }{\tag{$\Psc_{0}$}\label{opt-LQR}}{}
		\addConstraint{ \b x(0)}{=\b x_0}
		\addConstraint{ \b x'(t)}{=\b A(t) \b x(t) + \b B(t) \b u(t), \, \text{a.e.\ in} \, [0,T]}
		\addConstraint{ \b c_i(t)^\top \b x (t)}{\leq d_i(t),\,\forall \, t\in[0,T], \forall \,i\in\{1,\dots,P\},}
	\end{mini}
	where the state $\b x(t) \in \R^N$, the control $\b u(t) \in \R^M$, $\b A(t)\in \R^{N,N}$, $\b B(t)\in \R^{N,M}$, $\b C(t)=[\b c_1(t)^\top;\dots; \b c_P(t)^\top]\in \R^{P,N}$ ($\b c_i(t)\in\RR^N$), $\b d(t)=(d_i(t))_i\in \R^{P}$, while $\b Q(t)\in \R^{N,N}$ and $\b R(t) \in \R^{M,M}$ are positive semidefinite matrices.\\
	
	Below, for $q\in\{1, 2, \infty\}$, $L^q(0,T)$ denotes the $L^q$-space of functions over $[0,T]$ with integrable norms of the function values (resp.\ square-integrable, resp.\ bounded). We shall henceforth assume that, for all $t\in[0,T]$, $\b R(t) \succcurlyeq r \Id_M$ with $r>0$, as well as $\b A(\cdot)\in L^1(0,T)$, $\b B(\cdot)\in L^2(0,T)$, $\b Q(\cdot)\in L^1(0,T)$, and $\b R(\cdot)\in L^2(0,T)$. To have a finite objective, it is natural to restrict our attention to measurable controls satisfying $\b R(\cdot)^{1/2} \b u(\cdot)\in  L^2(0,T)$.\\ 
	
	Without state constraints, under mild assumptions, the unconstrained Linear Quadratic Regulator (LQR) enjoys an explicit solution defined through the Hamiltonian system and the related Riccati equation \citep[see e.g.][]{speyer2010primer}. With state constraints, little can be said as Pontryagin's Maximum Principle involves not only an adjoint vector but also measures supported on the constraint boundary. A comprehensive review of this approach can be found in \citet{hartl_survey_1995}. One has to guess beforehand when the state-constraint is active (at the so-called junction times) in order to write the first-order necessary condition \citep{hermant2009stab}. Secondly one has to impose assumptions to derive the magnitude of the discontinuities of the adjoint vector. This has proven to be intractable and made state-constrained continuous-time optimization a difficult problem. Let us provide an intuition for the appearance of discontinuities. If one follows an optimal trajectory of the LQR starting in the interior of state constraints, one may reach the boundary while the unconstrained Hamiltonian system of the Maximum Principle may incite to use a control leading to violation of the constraint. One has then to apply a different control to remain in the constraint set, possibly generating a discontinuity in the adjoint vector.
	
	Although LQR problems stand at the origin of control theory, research is still active in the field, not only for its numerous applications (see e.g.\ the examples of \citet{burachik2014duality} and references within), but also for its theoretical aspects, even without constraints \citep{bourdin2017linearquadratic} or just control constraints \citep{burachik2014duality}. Many of these improvements are motivated by model predictive control, considered for instance in a time-invariant discrete-time state-constrained setting in \citet{grune2018MPC} or continuous-time \citep{Keulen2020MPC}. In particular, \citet{kojima2004lq} proved that the solutions of a time-invariant LQR with discretized constraints converge to the solution of \eqref{opt-LQR}, putting emphasis on function spaces of controls. As a matter of fact, the aforementioned approaches \rv{focus on the control, used to obtain the trajectories}. In the present study, trajectories are instead at the core of the analysis.
	
	When seeking a continuous-time numerical solution, one has to face an infinite number of pointwise constraints, and has to either relax the \rv{computationally intractable} optimization problem or tighten it. \emph{Relaxing} means either enforcing the constraint only at a finite number of points, without guarantees elsewhere \citep[as with any time discretization method, e.g.\ ][]{kojima2004lq}, or through soft penalties \citep{gerdts2012virtual}, such as approximations of barrier functions \citep{dower2019barriers}.	\emph{Tightening} usually implies either choosing $\b u(\cdot)$ in a convenient subspace of $L^2(0,T)$, for instance the one of piecewise constant functions\footnote{This is known as \emph{sampled-data} or \emph{digital control} \citep{Ackermann1985}, the \emph{sampled-data} terminology does not refer to machine learning techniques.} or of splines with prescribed knots \citep{mercy2016MotionPlan}, or through hard penalties, such as logarithmic barriers \citep{chaplais2011design}. Let us illustrate the difference between relaxing and tightening state constraints. Consider the problem of a traffic regulator whose aim is to enforce a speed limit over a highway. The drivers for their part want to go as far as possible in a given time. Deploying speed cameras ensures at best that the speed constraint is satisfied  locally (relaxing). However if a smaller maximum speed is imposed at the camera level (tightening), then the cars cannot accelerate enough to break the speed limit before reaching the next camera. In a nutshell, the kernel methods framework we advocate allows to compute both a threshold and the resulting trajectories.\\
	
	Kernel methods, being related to Green's functions, belong to a branch of functional analysis. Their history was already sketched by \citet{aronszajn50theory}, to whom we owe the modern formulation of the theory. The regain of interest thanks to support vector machines \citep[see e.g.][]{scholkopf02learning} has reinstated kernel methods as the most principled technique in machine learning. There have been many attempts since then at bridging kernel methods and control theory. The two are already related in the works of \citet{parzen1970statis,kailath1971RKHS}. More recently \citet{steinke2008kernels, marco2017design} have considered kernels for control systems, mainly to encode the input $\b u(\cdot)$ or for system identification purposes (see e.g.\ the reviews of \citet{pillonetto2014kernel, chiuso2019system}). \rv{Kernels have also been applied to  approximate the Koopman operator over observers of uncontrolled nonlinear systems, in connection with spectral analysis \citep{OWilliams2015,Fujii2019} or for given controls \citep{Sootla2018}. The kernel Hilbert spaces have also been used to define suitable domains for operators \citep{rosenfeld2019dynamic, giannakis2019reproducing}. In most cases, the kernel is taken off-the-shelf, as with Gaussian kernels in connection with Bayesian inference  \citep{singh2018kernelStab, bertalan2019learning}.}
		
	On the contrary, departing from the prevalent perspective of using kernel methods as nonlinear embeddings, this article rekindles with a long standing tradition of engineering kernels for specific uses. This view has been mainly supported by the statistics community, especially in connection with splines and Sobolev spaces \citep{wahba90spline, heckman2012theory}.\footnote{Drawing inspiration from linear control theory, \emph{control theoretic splines} were devised \citep{egerstedt2009control, fujioka2013control}, in particular for path-planning problems \citep{kano2018Bsplines}. Possibly unbeknownst to non-kernel users, kernel theory, sometimes known as \emph{abstract splines}, is the natural generalization of splines \citep[see e.g.][]{aubin2020ifac}.} For \eqref{opt-LQR}, we show below that the quadratic objective paired with the linear dynamics encode the relevant kernel, which defines the Hilbert space of controlled trajectories. As kernel methods deal with a special class of Hilbert spaces, they are natural to consider for linear systems or for linearizations of nonlinear systems. Nonetheless the interactions run deeper. For instance we prove below that the controllability Gramian is \rv{directly related to} matrix-valued kernels, and we recover the transversality condition merely through a representer theorem. \rvS{This approach was further extended in \citet{aubin2020Riccati} to the connexion between the Linear-Quadratic matrix-valued kernel defined below and the dual Riccati equation. Since the Riccati equation is often used for the online purpose of finding the control, \citet{aubin2020Riccati} also discusses how the kernel formalism effectively allows for optimal synthesis, favoring an offline trajectory-focused viewpoint.} 
	
	Our main result is to show that through the theory of kernel methods we can solve a novel strengthened version of \eqref{opt-LQR}, without changing the function spaces involved. This is achieved by introducing a finite number of second-order cone constraints, stronger than the infinite number of affine \rv{constraints}. We furthermore show that the solution of the strengthened problem can be made arbitrarily close to the solution of \eqref{opt-LQR}, and that it enjoys a finite representation. One can thus exactly solve the continuous-time problem using only finite dimensional convex optimization solvers. This computable trajectory, which is interior to the affine state constraints, can also foster intuitions on the behavior of the optimal solution of \eqref{opt-LQR}. The considered tightening of \eqref{opt-LQR} relies solely on the kernel formalism. It was first introduced in \citet{aubin2020hard_nips} and thoroughly extended in \citet{aubin2020hard_SDP}.
	
	In Section \ref{sec_ProbForm} we present our strengthened problem and show how the LQR problem can be expressed as a regression problem over a vector-valued reproducing kernel Hilbert space. Section \ref{sec_revisitingLQ} details the consequences of this framework and identifies the corresponding Linear-Quadratic kernel. Our main result on approximation guarantees is stated in Section \ref{sec_bounds}. Section \ref{sec_numerics} discusses the implementation and the numerical behavior of the strengthened constraints. The setting is also extended to intermediate or terminal equality constraints, as in path-planning problems. The annex pertains to conditions ensuring the existence of interior trajectories for \eqref{opt-LQR}.
	\section{Theoretical preliminaries and problem formulation}\label{sec_ProbForm}
	
	In this section, we present the tools from the theory of kernel methods that we shall apply. We then introduce our strengthened problem with second-order cone constraints.\\
	
	\noindent\tb{Notations:} We use the shorthand $\iv{1}{P}=\{1,\dots,P\}$. $\RR_+^N$ is the subset of $\R^N$ of elements with nonnegative components. $\BB_N$ denotes the closed unit ball of $\R^N$ for the Euclidean inner product, $\b 1_N$ the vector of all-ones. For a matrix $\b A\in \R^{N,N}$, we write by $\|\b A\|$ its operator norm. $\Id_N$ is the identity matrix of $\R^{N,N}$.  We chose not to explicit the output space for the function spaces to avoid cumbersome notations, as it can be always deduced from the context. The space of functions with continuous derivatives up to order $s$ is denoted by $\C^s(0,T)$. For a function $K(\cdot,\cdot)$ \rv{defined over a subset of $\R\times \R$}, $\partial_1K(\cdot,\cdot)$ denotes the partial derivative w.r.t.\ the first variable. For a Hilbert space $(\Hk,\left<\cdot,\cdot\right>_{K})$, $\BB_K$ is the closed unit ball of $\Hk$, $\|\cdot\|_K$ denoting the corresponding norm. Given a subspace $V\subset\Hk$, we denote by $V^\perp$ its orthogonal complement w.r.t.\ $\left<\cdot,\cdot\right>_{K}$. A $\mu$-strongly convex function $\Lcal:\Hk\mapsto \R$ is a function satisfying, for all $f_1,f_2\in\Hk$, $\alpha\in[0,1]$, $ \Lcal(\alpha f_1+(1-\alpha)f_2) + \alpha(1-\alpha)\frac{\mu}{2}\|f_1-f_2\|_K^2\le \alpha \Lcal(f_1)  +(1-\alpha) \Lcal(f_2).$ 
	\begin{Definition}\label{def_vRKHS}
		Let $\T$ be a non-empty set. A Hilbert space $(\Hk,\left<\cdot,\cdot\right>_{K})$ of $\R^N$-vector-valued functions defined on $\T$ is called a vector-valued reproducing kernel Hilbert space (vRKHS) if there exists a matrix-valued kernel $K:\T \times \T \rightarrow \R^{N,N}$ such that the \emph{reproducing property} holds: for all $t \in \T,\, \b p\in\R^N $,  $K(\cdot,t)\b p \in \Hk$ and for all $\b f \in \Hk$, $\b p^{\top}\b f(t) = \left<\b f,K(\cdot,t)\b p\right>_{K}$.
	\end{Definition}
	Many properties of real-valued RKHSs have been known since \citet{aronszajn50theory}, the general theory having been developed by \citet{Schwartz1964}. By Riesz's theorem, an equivalent definition of a vRKHS is that, for every $t \in \T$ and $\b p\in\R^N$, the evaluation functional $\b f\in\Hk \mapsto \b p^\top \b f(t)\in\RR$ is continuous. There is also a one-to-one correspondence between the kernel $K$ and the vRKHS $(\Hk,\left<\cdot,\cdot\right>_{K})$ (see e.g.\ \citep[Theorem 2.6]{micheli_matrix-valued_2014}), hence modifying $\T$ or changing the inner product changes the kernel. We shall use several classical properties: by symmetry of the scalar product, the matrix-valued kernel has a Hermitian symmetry, i.e.\ $K(s,t)=K(t,s)^\top$ for any $s,t\in\T$. Moreover, if the vRKHS can be written as $\Hk=\H_0\oplus\H_1$, then $\H_0$ and $\H_1$ equipped with $\left<\cdot,\cdot\right>_{K}$ are also vRKHSs, as closed subspaces of $\Hk$ for $\|\cdot\|_{K}$, and their kernels $K_0$ and $K_1$ satisfy $K=K_0+K_1$.\\
	
	\noindent Let us define our candidate for a vRKHS, the space $\Sx$ of trajectories satisfying the dynamical system of \eqref{opt-LQR}:
	\begin{equation}
	\Sx:=\{\b x(\cdot)\,|\, \exists \, \b u(\cdot) \text{ s.t.\ } \b x'(t)=\b A(t) \b x(t) + \b B(t) \b u(t) \text{ a.e.\ and $\int_{0}^{T}\b u(t)^{\top} \b R(t) \b u(t) \d t <\infty$} \}.\label{def_Sx}
	\end{equation}
	 There is not necessarily a unique choice of $\b u(\cdot)$ for a given $\b x(\cdot) \in \Sx$ (for instance if $\b B(t)$ is not injective for some $t$). Therefore, with each $\b x(\cdot) \in \Sx$, we associate the control $\b u(\cdot)$ having minimal norm based on the pseudoinverse of $\b B(t)^{\ominus}$ of $\b B(t)$ for the $\R^M$-norm $\|\cdot\|_{\b R(t)}:=\|\b R(t)^{1/2}\cdot\|$:
	\begin{align}
	\b u(t) = \b B(t)^{\ominus}[\b x'(t) - \b A(t) \b x(t)] \, \text{ a.e.\ in} \,[0,T]. \label{def_u_as_X-X'} 
	\end{align}
	The vector space $\Sx$ has then a natural scalar-product. As a matter of fact, the expression
	\begin{equation}
	\left<\b x_1(\cdot),\b x_2(\cdot)\right>_{K}:= \b x_1(0)^\top \b x_2(0) + \int_{0}^{T}[\b x_1(t)^\top \b Q(t)\b x_2(t) + \b u_1(t)^\top \b R(t)\b u_2(t)] \d t 
		\label{def_K-norm}
	\end{equation}
	is bilinear and symmetric over $\Sx \times \Sx$.\footnote{The description by $\Sx$ of the optimization variables effectively pushes controls in the background while bringing forth trajectories as the main object of study. This describes \eqref{opt-LQR} more as a regression problem over $\Sx$ than as an optimal control problem over controls.} It is positive definite over $\Sx$ as $\|\b x(\cdot)\|^2_K = \b 0$ implies that $\b u(\cdot)\stackrel{a.e}{\equiv} 0$ and $\b x(0)= \b 0$, hence that $\b x(\cdot)\equiv \b 0$.	Combining \eqref{def_u_as_X-X'} and \eqref{def_K-norm}, we can express $\|\cdot\|_K$ as a Sobolev-like norm split into two semi-norms $\|\cdot\|_{K_0}$ and $\|\cdot\|_{K_1}$
	\begin{equation}
	\|\b x(\cdot)\|^2_K=\underbrace{\|\b x(0)\|^2 }_{\|\b x(\cdot)\|^2_{K_0}} + \underbrace{\int_{0}^{T}[\|\b x(t)\|^2_{\b Q(t)} + \|\b B(t)^{\ominus}(\b x'(t) - \b A(t) \b x(t))\|^2_{\b R(t)}] \d t}_{\|\b x(\cdot)\|^2_{K_1}}.\label{def_1k-norm}
	\end{equation}
	 By Cauchy-Lipschitz's theorem, $\|\b x(\cdot)\|_{K_0}=\|\b x(0)\|$ defines a norm over the finite-dimensional subspace $\Sx_0$ of trajectories with null quadratic cost (hence null control):
	\begin{equation}
	\Sx_0:=\{\b x(\cdot)\,|\, \int_{0}^{T}\b x(t)^{\top}\b Q(t) \b x(t) \d t=0 \text{ and } \b x'(t)=\b A(t) \b x(t), \,\text{ a.e.\ in } [0,T]  \}.\label{def_Sx0}
	\end{equation}
	Let us define its (infinite-dimensional) orthogonal complement $\Sx_u:=(\Sx_0)^{\perp}$ in $\Sx$ w.r.t.\ $\|\cdot\|_K$. From now on we equip $\Sx$ (resp.\ $\Sx_0$, $\Sx_u$) with $\|\cdot\|_K$ (resp.\ $\|\cdot\|_{K_0}$, $\|\cdot\|_{K_1}$). These will be shown to all be vRKHSs. Suppose we identified the matrix-valued kernels spawning them (a procedure to be found in Section \ref{sec_revisitingLQ}). In a recent approach in kernel methods \citep{aubin2020hard_nips}, developed for regression problems with constraints over derivatives, we suggested replacing ``$\b C(t) \b x (t)\leq \b d(t)\,(\forall \, t\in[0,T])$`` by the following strengthened second-order cone (SOC)\footnote{The "second-order cone" terminology is classical in optimization following the similarity between \eqref{def_SOC_cons} and the definition of the Lorentz cone $\{(\b z,r)\in\R^{N+1}\,|\, \|\b z\|\le r\}$.} constraints:
	\begin{align}\label{def_SOC_cons}
	\eta_i(\delta_m,t_m) \|\b x(\cdot)\|_K + \b c_i (t_m)^{\top} \b x (t_m) \leq d_{i}(\delta_m,t_m),\,\forall \, m\in \iv{1}{N_P}, \forall \,i\in\iv{1}{P}
	\end{align}
	where the $(t_m)_{m\in \iv{1}{N_P}}\in[0,T]^{N_P}$ are $N_P$ time points associated to radii $\delta_m>0$ satisfying $[0,T]\subset \cup_{m\in \iv{1}{N_P}} [t_m-\delta_m,t_m+\delta_m]$. The constants $\eta_i(\delta_m,t_m)$ and $d_{i}(\delta_m,t_m)$ are then defined as:
	\begin{align*}
	\eta_i(\delta_m,t_m) &:=\sup_{t\,\in\,[t_m-\delta_m,t_m+\delta_m]\cap [0,T]} \|K(\cdot,t_m)\b c_{i}(t_m)-K(\cdot,t)\b c_{i}(t) \|_K, \\
	d_{i}(\delta_m,t_m)&:=\inf_{t\,\in\,[t_m-\delta_m,t_m+\delta_m]\cap [0,T]} d_{i}(t).
	\end{align*}
	\rvS{This tightening of the constraints stems from interpreting $\eta_i(\delta_m,t_m) \|\b x(\cdot)\|_K$ as an upper bound of the modulus of continuity of the unknown $\b C(\cdot)\b x(\cdot)$ defined as follows
	\begin{align}\label{eq:mod_continuity_constraint}
	\omega^{(\b C\b x)}_i(\delta_m,t_m) &:=\sup_{t\,\in\,[t_m-\delta_m,t_m+\delta_m]\cap [0,T]} \underbrace{| \b c_i (t)^{\top} \b x (t)-  \b c_i (t_m)^{\top} \b x (t_m)|}_{|\langle x(\cdot),  K(\cdot,t)c(t)- K(\cdot,t_m)c(t_m)\rangle_K|}\le \eta_i(\delta_m,t_m) \|\b x(\cdot)\|_K.
	\end{align}
	This inequality is obtained applying successively the reproducing property and the Cauchy–Schwarz inequality. Since the intractable modulus of continuity controls the variations of $\b C(t) \b x (t)$, the SOC upper bound provides a tractable tightening. This interpretation through moduli of continuity was extensively developed in Section 3.2 of \citet{aubin2020hard_SDP}.}

	With the above notations, our strengthened time-varying linear quadratic optimal control problem with finite horizon and finite number of SOC constraints is\footnote{Since $\|\b x(\cdot)\|^2_{K_0}=\|\b x(0)\|^2$ and $\b x(0)$ is fixed, we replaced the integral $\|\b x(\cdot)\|^2_{K_1}$ by $\|\b x(\cdot)\|^2_{K}$ in the objective. Recently, in \citet{aubin2020Riccati}, another choice of inner product was introduced, by incorporating  quadratic terminal costs $g$ into $\|\b x(\cdot)\|^2_{K}$ for another choice of inner product than \eqref{def_K-norm}.}
	\begin{mini}
		{\substack{\b x(\cdot)\in\Sx,\\ \b x(0) =\b x_0}}{\mathrlap{g(\b x(T))+\|\b x(\cdot)\|^2_{K}} \phantom{\eta_i(\delta_m,t_m) \|\b x(\cdot)\|_K + \b c_i (t_m)^{\top} \b x (t_m) \leq d_{i}(\delta_m,t_m),\, \forall \, m\in \iv{1}{N_P}, \forall \,i\in\iv{1}{P}}}
		{\label{opt-LQR_SOC}\tag{$\Psc_{\text{$\delta$,fin}}$}}{}
		\addConstraint{\eta_i(\delta_m,t_m) \|\b x(\cdot)\|_K + \b c_i (t_m)^{\top} \b x (t_m)}{\leq d_{i}(\delta_m,t_m),\, \forall \, m\in \iv{1}{N_P}, \forall \,i\in\iv{1}{P}.}
	\end{mini}
	The introduction of \eqref{opt-LQR_SOC} as an approximation of \eqref{opt-LQR} entirely stems from the vRKHS formalism and does not result from optimal control considerations. It relies on an inner approximation of a convex set in an infinite-dimensional Hilbert space \citep[see][Section 3.1]{aubin2020hard_SDP}. \rv{From a machine learning perspective, the initial condition and the terminal cost act as a ``loss`` function, whereas the quadratic cost is turned into a norm over $\Sx$ and can thus be interpreted as a ``regularizer``. Departing also from optimal control, the tightening is obtained by incorporating the quadratic part of the objective \eqref{def_1k-norm} in the state constraints to form \eqref{def_SOC_cons}. As discussed in \citet{aubin2020hard_nips}, introducing \eqref{def_SOC_cons} leads to a finite number of evaluations of the variable $\b x(\cdot)$ in \eqref{opt-LQR_SOC} which allows for a representer theorem (Theorem \ref{thm_representer} in Section \ref{sec_Q0}).}\\
	
	Our goal is to show that \eqref{opt-LQR_SOC} is indeed a tightening of \eqref{opt-LQR}, enjoying a representer theorem providing a finite-dimensional representation of the solution of problem \eqref{opt-LQR_SOC}. We also bound the distance between the trajectories solutions of \eqref{opt-LQR} and \eqref{opt-LQR_SOC}, and prove that it can be made as small as desired by refining the time-discretization grid.
	
	\section{Revisiting LQ control through kernel methods}\label{sec_revisitingLQ}
	
	This section \rv{presents a step-by-step approach to identify} the matrix-valued kernel $K$ of the Hilbert space $\Sx$ of solutions of a \rv{linear} control system, equipped with the scalar product \eqref{def_K-norm}. This is done independently from the state constraints which effect is only to select a closed convex subset of the space of trajectories. In Section \ref{sec_Q0}, we consider the case $\b Q \equiv \b 0$ which enjoys explicit formulas. \rv{We also express a representer theorem (Theorem \ref{thm_representer} suited for problems of the form \eqref{opt-LQR_SOC}.} This allows us to revisit, through the kernel framework, classical notions, such as the solution of the unconstrained LQR problem, or the definition of the Gramian of controllability. In Section \ref{sec_Qnot0}, we consider the case $\b Q \not\equiv \b 0$ and relate our solution to an adjoint equation over matrices. Furthermore, the identification of kernels developed \rv{in Section \ref{sec_revisitingLQ}} is by no means restricted to finite $T$, hence the kernel formalism can also tackle infinite-horizon problems.\footnote{\rv{We need to assume $T$ to be finite for the results of Section \ref{sec_bounds} to hold, as we use a procedure based on compact coverings to deal with the state constraints.}} 
	
	Let us denote by $\b \Phi_{\b A}(t,s)\in\R^{N,N}$ the state-transition matrix of $\b x'(\tau)=\b A(\tau) \b x(\tau)$, defined from $s$ to $t$. The key property used throughout this section is \rv{the variation of constants \rvS{formula}, a.k.a.\ Duhamel's principle,} stating that for any absolutely continuous $\b x(\cdot)$ such that $ \b x'(t)=\b A(t) \b x(t) + \b B(t) \b u(t)$ a.e.\, we have
	\begin{equation}
	\b x(t) = \b \Phi_{\b A}(t,0) \b x(0)+\int_{0}^{t} \b \Phi_{\b A}(t,\tau) \b B(\tau) \b u(\tau) \d \tau.\label{eq_Duhamel}
	\end{equation}
	
	\begin{Lemma}\label{lem_SxHilbert} $(\Sx,\left<\cdot,\cdot\right>_{K})$ is a vRKHS.
	\end{Lemma}		
	\noindent\textbf{Proof:} We have to show that: $i)$ $(\Sx,\left<\cdot,\cdot\right>_{K})$ is a Hilbert space, $ii)$ for every $t \in [0,T]$ and $\b p\in\R^N$, the evaluation functional $\b x(\cdot)\in\Sx \mapsto \b p^\top \b x(t)\in\RR$ is continuous.
	
	\noindent  $i)$ From \eqref{def_K-norm} and the discussion of Section \ref{sec_ProbForm}, it is obvious that $\left<\cdot,\cdot\right>_{K}$ is a scalar product. We just have to show that $\Sx$ is complete. Let $(\b x_n(\cdot))_n$ be a
	Cauchy sequence \rv{in $\Sx$}, with associated controls $(\b u_n(\cdot))_n$. Then $(\|\b x_n(\cdot)\|_K)_n$ \rv{is a Cauchy sequence in $\R$ and thus converges}, so $(\|\b x_n (0)\|)_n$ and $(\|\b R(\cdot)^{1/2} \b u_n(\cdot)\|_{L^2(0,T)})_n$ are bounded. Since $\b R(t) \succcurlyeq r \Id_M$ with $r>0$, $(\| \b u_n(\cdot)\|_{L^2(0,T)})_n$ is thus bounded too, and we can take a subsequence $(\b u_{n_i}(\cdot))_i$ weakly converging to some $\b u(\cdot)$. Let $s,t\in[0,T]$, 
	\begin{align}
	\b x_n(t) - \b x_n(s) \stackrel{\eqref{eq_Duhamel}}{=} (\b \Phi_{\b A}(t,s) - \Id_N)\b x_n(s)+\int_{s}^{t} \b \Phi_{\b A}(t,\tau) \b B(\tau) \b u_n(\tau)  \d \tau.\label{eq_Duhamel_diff}
	\end{align}
	Taking $s=0$, as $(\|\b x_n (0)\|)_n$ is bounded, $\b A(\cdot)\in L^1(0,T)$ and $\b B(\cdot)\in L^2(0,T)$, we have that $\b \Phi_{\b A}(\cdot,\cdot)$ is continuous and $\{\b x_n(\cdot
	)\}_n$ is uniformly bounded in $\C(0,T)$. Thus \eqref{eq_Duhamel_diff} implies that the sequence $(\b x_{n_i}(\cdot))_i$ is equicontinuous. By Ascoli's theorem, we can take a $(\b x_{n_{i_j}}(\cdot))_j$ uniformly converging to some $\b x(\cdot)$ satisfying \eqref{eq_Duhamel} for $\b u(\cdot)$, thus $\b x(\cdot)\in\Sx$.
	
	\noindent  $ii)$ Let $t \in [0,T]$, $\b p\in\R^N$, and $\b x(\cdot)\in\Sx$. By \eqref{eq_Duhamel} and Cauchy-Schwarz inequality,
	\begin{align*}
	\b p^\top \b x(t) &= \b p^\top \b \Phi_{\b A}(t,0) \b x(0)+\int_{0}^{t} \b p^\top \b \Phi_{\b A}(t,\tau) \b B(\tau) \b u(\tau) \d \tau\\
	|\b p^\top \b x(t)| &\le \sup_{\tau\in[0,T]} \| \b p^\top \b \Phi_{\b A}(t,\tau) \| \cdot\left(\| \b x(0)\|+\|\int_{0}^{t} \b B(\tau) \b u(\tau) \d \tau\|\right)\\
	&\leq \sup_{\tau\in[0,T]} \sqrt{2}\| \b p^\top \b \Phi_{\b A}(t,\tau) \| \cdot\left(\| \b x(0)\|^2+ \| \b B(\cdot)\|_{L^2(0,T)}^2 \| \b u(\cdot)\|_{L^2(0,T)}^2 \right)^{\frac{1}{2}}\\
	&\leq \sup_{\tau\in[0,T]} \sqrt{2}\| \b p^\top \b \Phi_{\b A}(t,\tau) \| \cdot\left(1+\frac{\| \b B(\cdot)\|_{L^2(0,T)}}{\sqrt{r}}\right)  \left(\| \b x(0)\|^2+ \| \b R(\cdot)^{\frac{1}{2}}\b u(\cdot)\|_{L^2(0,T)}^2 \right)^{\frac{1}{2}}.
	\end{align*}
	Hence \rv{the linear map} $\b x(\cdot)\in\Sx \mapsto \b p^\top \b x(t)\in\RR$ is continuous. 
	\begin{flushright}
		$\blacksquare$
	\end{flushright}
	By Definition \ref{def_vRKHS}, we know that a matrix-valued reproducing kernel $K(\cdot,\cdot)$ exists. We now repeatedly use \eqref{eq_Duhamel} to identify it.
	
	\subsection{Case $\b Q\equiv\b 0$}\label{sec_Q0}
	The space $\Sx$ of controlled trajectories is defined as in \eqref{def_Sx} equipped with the quadratic norm
	\begin{equation}
	\|\b x(\cdot)\|^2_K=\|\b x(0)\|^2 + \int_{0}^{T}\b u(t)^{\top} \b R(t) \b u(t) \d t=\| \b x(\cdot)\|_{K_0}^2+\| \b x(\cdot)\|_{K_1}^2, \label{def_1k-norm_Q0}
	\end{equation}
	where $\b u(\cdot)$ is defined as in \eqref{def_u_as_X-X'}. We can further \rv{specify} its subspaces	
	\begin{align*}
	\Sx_0=\{\b x(\cdot)\,|\, \b x'(t)=\b A(t) \b x(t), \,\text{a.e.\ in } [0,T]  \} \quad \quad 
	\Sx_u=\{\b x(\cdot)\,|\, \b x(\cdot)\in\Sx \text{ and } \b x(0)=0 \}. 
	\end{align*}
	By uniqueness of the reproducing kernel, we only have to exhibit a candidate $K(s,t)$ satisfying $\b p^{\top}\b x(t) = \left<\b x(\cdot),K(\cdot,t)\b p\right>_{K}$ ($\forall t \in [0,T], \b x \in \Sx$, $\b p\in\R^N$) and $K(\cdot,t)\b p \in \Sx$ ($\forall t \in [0,T]$, $\b p\in\R^N$). The space $(\Sx_0, \|\cdot\|_{K_0})$ being finite dimensional, we can right away identify its kernel\footnote{This is a classical result for finite dimensional vRKHSs. Fix any family $\{v_j\}_{j\in  \iv{1}{N}}$ spanning $\Sx_0$, let $\b V(s):= [v_j(s)]_{j\in  \iv{1}{N}}\in\R^{N,N}$ and $\b G_v := (\langle v_i,v_j \rangle_{K_0})_{i,j\in  \iv{1}{N}}$. The matrix $\b G_v$ is invertible as $\|\cdot\|_{K_0}$ is a norm over $\Sx_0$ and thus $K_0(s,t)=\b V(s)^{\top} \b G_v^{-1} \b V(t)$. Here we have $\b V(s)=\b \Phi_{\b A}(s,0)^\top$ and $\b G_v =\Id_N$.}
	\begin{equation}
	K_0(s,t)=\b \Phi_{\b A}(s,0) \b \Phi_{\b A}(t,0)^{\top}.\label{def_K0}
	\end{equation}
	As $\Sx=\Sx_0 \oplus \Sx_u$, from the properties of sums of kernels, we derive that we should look for $K$ of the form $K_0+K_1$ for which the reproducing property, with  $\left<\cdot,\cdot\right>_{K}$ defined in \eqref{def_1k-norm}, writes \rv{as follows, for all $t \in [0,T]$, $\b p\in\R^N$, $\b x(\cdot)\in\Sx$,}
	\begin{equation}
	\b p^{\top}\b x(t)=(K(0,t)\b p)^{\top} \b x(0) + \int_{0}^{T}  \left[(\b B(t)^{\ominus}(\b \partial_1 (K(s,t)\b p )- \b A(t) K(s,t)\b p))^{\top} \b R(s) \b u(s) \right] \d s. \label{eq_repro_Q0}
	\end{equation}
	 Setting $\partial_1 K(s,t):\b p\mapsto\partial_1 (K(s,t)\b p)$, let us define formally $\b U_t(s) := \b B(s)^{\ominus}(\b \partial_1 K_1(s,t) - \b A(s) K_1(s,t))$. By the Hermitian symmetry of $K$ and the fact that $K_0(\cdot,t)\b p$ belongs to $\Sx_0$ and $K_1(0,t)\b p=\b 0$, \rv{\eqref{eq_repro_Q0} holds if and only if for all $t \in [0,T]$ and $\b x(\cdot)\in\Sx$}
	\begin{equation*}
	\b x(t)=K_0(t,0) \b x(0) + \int_{0}^{T}  \b U_t(s)^{\top} \b R(s) \b u(s) \d s. 
	\end{equation*}
	This expression can be identified with \eqref{eq_Duhamel} when defining $\b U_t(s)$ as follows
	\begin{equation}\label{def_Uts}
	\b U_t(s):=	\left\{\begin{array}{cl}
		\b R(s)^{-1}\b B(s)^{\top} \b \Phi_{\b A}(t,s)^{\top}\; &\forall s\leq t,\\ 
		\b 0\; &\forall s > t.
		\end{array} \right.
	\end{equation}
	Consequently, to ensure that, for all $t \in [0,T]$, $\b p\in\R^N$,  $K_1(\cdot,t)\b p \in \Sx_u$, $K_1$ has to satisfy the following differential equation for any given $t\in[0,T]$:
	\begin{equation}
	\b \partial_1 K_1(s,t) = \b A(s) K_1(s,t) + \b B(s)\b U_t(s) \text{ a.e.\ in } [0,T] \text{ with } K_1(0,t)=\b 0.\label{def_EDO_K1}
	\end{equation}
	Since $\b A(\cdot)\in L^1(0,T)$, $\b B(\cdot)\in L^2(0,T)$, and $\b R(t) \succcurlyeq r \Id_M$ with $r>0$, we have that $\b U_t(\cdot)\in L^2(0,T)$. By applying \rv{the variation of constants \rvS{formula}} to \eqref{def_EDO_K1}, with $\b U_t(\cdot)$ defined in \eqref{def_Uts}, we get an explicit expression for $K_1$, satisfying a Hermitian symmetry when permuting $s$ and $t$,
	\begin{equation}
	K_1(s,t) = \int_{0}^{\min(s,t)} \b \Phi_{\b A}(s,\tau) \b B(\tau) \b R(\tau)^{-1} \b B(\tau)^{\top} \b \Phi_{\b A}(t,\tau)^{\top}  \d \tau.\label{def_K1}
	\end{equation}
	\textbf{Remark (Gramian):} Formula \eqref{def_K1} for $K_1(T,T)$ corresponds to the Gramian of controllability. The link is straightforward as the controllability problem of steering a point from $(0,\b 0)$ to $(T,\b x_T)$ simply writes as, with $\b u(\cdot)$ defined as in \eqref{def_u_as_X-X'} and $\b R(\cdot) \equiv \Id_M$,
		\begin{mini}
		{\substack{\b x(\cdot)\in\Sx}}{\int_{0}^{T} \|\b u(t)\|^2 \d t}{\label{opt-LQ_Gram}}{}
		\addConstraint{ \b x(0)}{=\b 0}
		\addConstraint{ \b x(T)}{=\b x_T},
		\end{mini}
	which set of solutions can actually be made explicit.\footnote{When choosing \rv{the terminal cost} $g(\cdot)$ \rv{to be} equal to the indicator function of $\{\b x_T\}$, \eqref{opt-LQ_Gram} does correspond to \eqref{opt-LQR} in the absence of state constraints.} As a matter of fact, in kernel methods, it is classical to look for a ``representer theorem``, i.e.\ a necessary condition to ensure that the solutions of an optimization problem live in a finite dimensional subspace of $\Sx$ and consequently enjoy a finite representation. Such theorems are usually stated without constraints and for real-valued kernels \citep[e.g.][]{Scholkopf2001}. Here we formulate a representer theorem for conic constraints and matrix-valued kernels, as it will prove instrumental to derive a finite formulation for the SOC-strengthening of \eqref{opt-LQR}.
	\begin{Theorem}[Representer theorem] \label{thm_representer} Let $(\Hk,\left<\cdot,\cdot\right>_{K})$ be a vRKHS defined on a set $\T$. Let $P\in\NN$ and, for $i\in\iv{0}{P}$ and given $N_i\in\NN$, $\{t_{i,j}\}_{j\in \iv{1}{N_i}}\subset \T$. Consider the following optimization problem with ``loss`` function $L:\RR^{N_0}\rightarrow \RR\cup\{+\infty\}$, strictly increasing ``regularizer`` function $\Omega:\RR_+\rightarrow \RR$, and constraints $d_i:\RR^{N_i}\rightarrow \RR$, $\lambda_i\geq 0$ and $\{\b c_{i,m}\}_{m\in\iv{1}{N_i}} \subset\R^N$,
		\begin{argmini*}
			{\substack{\b f\in\Hk}}{ L\left(\b c_{0,1}^\top \b f(t_{0,1}),\dots, \b c_{0,N_0}^\top \b f(t_{0,N_0})\right) + \Omega\left(\|\b f \|_K\right)}{}{\bar{\b f}\in}
			\addConstraint{ \lambda_i\|\b f \|_K}{\leq d_i(\b c_{i,1}^\top \b f(t_{i,1}),\dots, \b c_{i,N_i}^\top\b f(t_{i,N_i})),\, \forall\, i\in\iv{1}{P}.}
		\end{argmini*}
	Then, for any minimizer $\bar{\b f}$, there exists $\{\b p_{i,m}\}_{m\in \iv{1}{N_i}}\subset \RR^N$ such that $\bar{\b f}= \sum_{i=0}^{P} \sum_{m=1}^{N_i} K(\cdot,t_{i,m}) \b p_{i,m}$ with $\b p_{i,m}=\alpha_{i,m} \b c_{i,m}$ for some $\alpha_{i,m}\in\RR$.
	\end{Theorem}
	\noindent\textbf{Proof:} Let $\bar{\b f}$ be an optimal solution and let $V := \Sp\left(\{K(\cdot,t_{i,m})  \b c_{i,m}\}_{m\le N_i, \, i\le P}\right)$.
	Take $\b v\in V$ and $\b w\in V^\perp$ such that $\bar{\b f}=\b v+\b w$. As $ \b c_{i,m}^\top \b w(t_{i,m})=\langle\b w(\cdot), K(\cdot,t_{i,m})  \b c_{i,m}\rangle_K=0$, the terms appearing in $L$  and $d_i$ are the same for  $\bar{\b f}$ and $\b v$. Moreover, $\|\b v \|_K \le \|\bar{\b f} \|_K$, hence $\b v$ belongs to the constraint set since $\bar{\b f}$ does. Furthermore $\Omega\left(\|\b v \|_K\right)\le \Omega\left(\|\bar{\b f} \|_K\right)$, so, by optimality of $\bar{\b f}$, $\b w=0$ which concludes the proof.
\begin{flushright}
	$\blacksquare$
\end{flushright}
	In other words, Theorem \ref{thm_representer} states that to each time $t$ where the variable $\b f$ is evaluated corresponds a multiplier $\b p_t \in \R^N$ in the expression of the optimal solutions. Hence, if the number of such evaluations is finite, then the representation of $\bar{\b f}$ is finite.\footnote{The property of having a finite number of evaluations is precisely what distinguishes the unconstrained controllability problem \eqref{opt-LQ_Gram} or the SOC-constrained problem \eqref{opt-LQR_SOC} from the original state-constrained problem \eqref{opt-LQR} which has an infinite number of affine constraints.} Besides, a representer theorem, like Pontryagin's Maximum Principle, is only a necessary condition on the form of the solutions. Theorem \ref{thm_representer} guarantees the existence or uniqueness of an optimal solution only when coupled with other assumptions (e.g.\ that $L$ and $\Omega$ are convex, and $L$ is lower semi-continuous). This further \rv{highlights} the analogy with the Maximum Principle in the quadratic case. Theorem \ref{thm_representer} is instrumental to obtain a finite-dimensional equivalent of \eqref{opt-LQR_SOC}.\\
	
	Theorem \ref{thm_representer} applied to \eqref{opt-LQ_Gram} implies that any candidate optimal solution $\bar{\b x}(\cdot)$ can be written as $\bar{\b x}(s)=K(s,0) \b p_0+K(s,T) \b p_T$, with $ \b p_0,  \b p_T\in\R^N$. As $\bar{\b x}(0)=\b 0$, $\text{proj}_{\Sx_0}(\bar{\b x}(\cdot))=\b 0$, and as $K_1(\cdot,0)\equiv \b 0$, $\bar{\b x}(s)=K_1(s,T) \b p_T$. So $\bar{\b x}(T)=\b x_T$ is satisfied if and only if $\b x_T\in\text{Im}(K_1(T,T))$ where the operator $K_1(T,T)$ is defined by \eqref{def_EDO_K1}, setting $\b R\equiv \text{Id}$. Hence \eqref{opt-LQ_Gram} has a solution for any $\b x_T$ (i.e.\ the system is controllable) if and only if the Gramian of controllability $K_1(T,T)$ is invertible.\\
	
	\noindent\textbf{Remark (LQR without state constraints):} We derive also from the kernel framework the transversality condition, as well as the classical solution of the LQR problem without state constraints, defined as follows, with $\b u(\cdot)$ again defined as in \eqref{def_u_as_X-X'},
	\begin{mini}
		{\substack{\b x(\cdot)\in\Sx}}{g(\b x(T))+\frac{1}{2}\int_{0}^{T} \b u(t)^{\top} \b R(t) \b u(t) \d t}{\tag{$\Psc_{\text{uncons}}$}\label{opt-LQ_nocons}}{}
		\addConstraint{ \b x(0)}{=\b 0.}
	\end{mini}
	Similarly, through the representer theorem, we deduce that $\bar{\b x}(\cdot)=K_1(\cdot,T) \b p_T$. Hence, by the reproducing property,
	$$\int_{0}^{T} \bar{\b u}^{\top}(t) \b R(t) \bar{\b u}(t) \d t=\|K_1(\cdot,T) \b p_T\|_K^2=\b p_T^{\top} K_1(T,T)\b p_T.$$  
	Assume that $g(\cdot)\in\mathcal{C}^1(\R^N,\R)$ \rv{and that it is convex}. Applying the first-order optimality condition, we conclude that
	\begin{align}
		\b 0&=\nabla\left(\b p\mapsto g(K_1(T,T)\b p)+\frac{1}{2} \b p^{\top} K_1(T,T)\b p \right)(\b p_T)=K_1(T,T)(\nabla g(K_1(T,T)\b p_T) +\b p_T).\label{eq_transversality}
	\end{align}
	So it is sufficient to take $\b p_T= -\nabla g(K_1(T,T)\b p_T) = -\nabla g(\bar{\b x}(T))$, i.e.\ to have the transversality condition satisfied. However this formula says more than that, as it covers the problem of degeneracies of the ``controllability Gramian`` $K_1(T,T)$ and gives an explicit equation \eqref{eq_transversality} to be satisfied by $\b p_T$. Notice that we do not consider any adjoint equation, only adjoint vectors that are not explicitly propagated. In our framework, the Hamiltonian is implicit.
	\subsection{Case $\b Q\not\equiv\b 0$}\label{sec_Qnot0}
	For the case with $\b Q \not\equiv \b 0$, we have a more intricate formula. The reproducing property for $K$, \rvS{in which the term $\b U_t(\cdot)$ will be explicitly specified below}\footnote{Recall that $\b U_t(s) := \b B(s)^{\ominus}[\b \partial_1 K(s,t) - \b A(s) K(s,t)]$.}, writes as follows, \rv{for all $t \in [0,T]$, $\b p\in\R^N$, $\b x(\cdot)\in\Sx$,}
	\begin{equation}
	\b p^{\top}\b x(t)=(K(0,t)\b p)^{\top} \b x(0) + \int_{0}^{T} (K(s,t)\b p)^{\top} \b Q(s) \b x(s) \d s+\int_{0}^{T} (\b U_t(s)\b p))^{\top} \b R(s) \b u(s) \d s. \label{eq_repro_Q}
	\end{equation}
	By the Hermitian symmetry of $K$ and \rv{the variation of constants \rvS{formula}} \eqref{eq_Duhamel}, we can rewrite \eqref{eq_repro_Q} as, \rv{for all $t \in [0,T]$, $\b x(\cdot)\in\Sx$,}
	\begin{align*}
	\b x(t)=K(t,0) \b x(0) &+  \int_{0}^{T} K(t,s) \b Q(s) \left(\b \Phi_{\b A}(s,0) \b x(0)+\int_{0}^{s} \b \Phi_{\b A}(s,\tau) \b B(\tau) \b u(\tau)  \d \tau\right) \d s\\
	 &+\int_{0}^{T}  \b U_t(s)^{\top} \b R(s) \b u(s) \d s.
	\end{align*}
	After some regrouping of terms, a change of integration bounds and identification with \eqref{eq_Duhamel}, we get the integral equations:
	\begin{align*}
	\b \Phi_{\b A}(t,0)&= K(t,0) +\int_0^{T} K(t,s)  \b Q(s) \b \Phi_{\b A}(s,0) \d s,\\
	\forall\, s\leq t,\; \b \Phi_{\b A}(t,s)\b B(s)&= \b U_t(s)^{\top} \b R(s) +\int_{s}^{T} K(t,\tau)  \b Q(\tau) \b \Phi_{\b A}(\tau,s) \b B(s)  \d \tau \\
	\forall\, s> t,\; \b 0&= \b U_t(s)^{\top} \b R(s) +\int_{s}^{T} K(t,\tau)  \b Q(\tau) \b \Phi_{\b A}(\tau,s) \b B(s)  \d \tau,
	\end{align*}
	which can be summarized as 
	\begin{align}
	K(t,0)=& \; \b \Phi_{\b A}(t,0)-\tilde{K}(t,0)\text{ with } \tilde{K}(t,s):=\int_{s}^{T} K(t,\tau)  \b Q(\tau) \b \Phi_{\b A}(\tau,s)  \d \tau, \nonumber \\
	\b U_t(s)^{\top} \b R(s)=&	\left\{\begin{array}{cl}
	(\b \Phi_{\b A}(t,s)-\tilde{K}(t,s))\b B(s)\; &\forall s\leq t,\\ 
	-\tilde{K}(t,s)\b B(s)\; &\forall s > t.
	\end{array} \right.\label{eq_integralK}
	\end{align}

	Although not as explicit as formula \eqref{def_K1} from the case $\b Q\equiv\b 0$, this integral expression for $\b Q\not\equiv\b 0$ will still prove valuable to investigate the regularity of $K(\cdot,\cdot)$ (Lemma \ref{lem_modCo_noyau} below). 
	To provide further insight on this expression, again for fixed $t$, let us introduce formally an adjoint equation \rv{for} a variable $\b \Pi(s,t)\in\RR^{N,N}$,
	\begin{align}\label{eq:EDO_adjoint}
		\partial_1\b \Pi(s,t)=- \b A(s)^\top \b \Pi(s,t) +\b  Q(s)K(s,t) \quad \quad \b  \Pi(T,t)=\Id_N.
	\end{align}
	Again, applying \rv{the variation of constants \rvS{formula}} to $\b \Pi(s,t)$, taking the transpose and owing to the symmetries of $K(\cdot,\cdot)$ and $\b \Phi(\cdot,\cdot)$, we derive that
	\begin{align*}
	\b \Pi(s,t) &= \b \Phi_{(-\b A^\top)}(s,T) \b \Pi(T,t)+\int_{T}^{s} \b \Phi_{(-\b A^\top)}(s,\tau) \b  Q(\tau)K(\tau,t)  \d \tau\\
	\b \Pi(s,t)^\top &= \b \Phi_{\b A}(T,s)-\int_{s}^{T}  K(t,\tau) \b  Q(\tau) \b \Phi_{\b A}(\tau,s)  \d \tau= \b \Phi_{\b A}(T,s)-\tilde{K}(t,s).
	\end{align*}
	Since $\partial_1 K(s,t) = \b A(s) K(s,t) + \b B(s) \b U_t(s)$, by \eqref{eq_integralK}, for any given time $t$,
	\begin{equation}\label{eq_Hamilton_full}
	\left.\begin{array}{rl}
	\partial_1 K(s,t) &= \b A(s) K(s,t) + \b B(s) \b R(s)^{-1} \b B(s)^\top \left\{
	\begin{array}{l}
	\b \Pi(s,t) - \b \Phi_{\b A}(T,s)^\top + \b \Phi_{\b A}(t,s)^\top \;\forall s\leq t,\\ 
	\b \Pi(s,t) - \b \Phi_{\b A}(T,s)^\top\quad\forall s > t.
	\end{array} \right.\\ 
	K(0,t)&=\b  \Pi(0,t) +  \b \Phi_{\b A}(t,0)^\top -  \b \Phi_{\b A}(T,0)^\top.  	\end{array}\right.
	\end{equation}

	The difference of behavior between the two cases $\b Q\equiv\b 0$ and $\b Q\not\equiv\b 0$ is classical in optimal control. While the control equation runs forward in time, the adjoint equation runs backward. For $\b Q\equiv\b 0$, the adjoint equation can be solved independently from $K$, which is why $\b \Pi_t$ was not introduced. For $\b Q \not \equiv\b 0$, we have two coupled differential equations \eqref{eq:EDO_adjoint}-\eqref{eq_Hamilton_full} over $K(\cdot,t)$ and $\b \Pi_t(\cdot)$. This system does not enjoy an explicit expression, however its solutions can still be computed as a two-point boundary value problem. For quadratic terminal costs, see \citet{aubin2020Riccati}, where the connections between computing the kernel and solving the Hamiltonian system or its Riccati equation counterpart are stressed.
	
	\section{Theoretical approximation guarantees}\label{sec_bounds}
	In this section, we show that the SOC-constrained problem  \eqref{opt-LQR_SOC} is a tightening of the original problem \eqref{opt-LQR}. We also provide bounds on the $\|\cdot\|_K$-distance between the optimal trajectory of \eqref{opt-LQR} and that of \eqref{opt-LQR_SOC}. This shows that the SOC-tightening is consistent in a numerical analysis sense, as, for bounded kernels $K$, convergence in $\|\cdot\|_K$ is stronger than uniform convergence of the states, and also implies convergence of the $L^2$-norms of the controls. We prove that the kernels $K(\cdot,\cdot)$ identified in Section \ref{sec_revisitingLQ} \rv{are} indeed $C^0$-continuous.\\
	
	We shall manipulate various forms of state constraints. We thus write our generic problem \eqref{opt-LQR_gen}, with \rv{objective $\Lcal$ and} constraints defined through the constraint set $\V_*$, as follows
	\begin{argmini}
		{\substack{\b x(\cdot)\in \V_*}}{\rv{\Lcal(\b x(\cdot)):=\;}g(\b x(T))+\|\b x(\cdot)\|^2_{K}}{\tag{$\Psc_{*}$}}{ \bar{\b x}^* (\cdot)\in\label{opt-LQR_gen}}
		\addConstraint{ \b x(0)}{=\b x_0.}
	\end{argmini}
	Existence and uniqueness of the solution for each $(\Psc_{*})$ will be discussed below. Recall that $\Sx$ is defined in \eqref{def_Sx}. Let $(t_m)_{m\in\iv{1}{N_0}}\in[0,T]^{N_0}$ be $N_0$ time points associated to radii $\delta_m>0$ chosen so that they form a covering $[0,T]\subset\cup_{m\in\iv{1}{N_0}} [t_m-\delta_m,t_m+\delta_m]$. The vectors $\b d_{m}(\delta_m,t_m):=(d_i(\delta_m,t_m))_{i\in\iv{1}{P}}\in \R^{P}$, $\etaArr(\delta,t):=(\eta_i(\delta,t))_{i\in\iv{1}{P}}\in \R^{P}$ and $\omeArr(\delta,t):=(\omega_i(\delta,t))_{i\in\iv{1}{P}}\in \R^{P}$ are defined component-wise:\footnote{The computation of $\eta_i$ can be performed using that, by the reproducing property, $\|K(\cdot,t)\b c_{i}(t)-K(\cdot,s)\b c_{i}(s) \|^2_K=\b c_{i}(t)^\top K(t,t)\b c_{i}(t)+\b c_{i}(s)^\top K(s,s)\b c_{i}(s)- 2\b c_{i}(t)^\top K(t,s) \b c_{i}(s)$. \rv{We chose to overload the notation of $d_i$ to define the constants $d_{i}(\delta_m,t_m)$ in order to draw the parallel with the other perturbations of the constraints, $\eta_i(\delta,t)$ and $\omega_i(\delta,t)$.} }
	\begin{align}
	\eta_i(\delta,t) &:=\sup_{s\,\in\,[t-\delta,t+\delta]\cap[0,T]} \|K(\cdot,t)\b c_{i}(t)-K(\cdot,s)\b c_{i}(s) \|_K,\label{def_eta}\\
	\omega_i(\delta,t) &:=\sup_{s\,\in\,[t-\delta,t+\delta]\cap[0,T]} |d_{i}(t)- d_{i}(s) |,\label{def_omega}\\
	d_{i}(\delta_m,t_m)&:=\inf_{s\,\in\,[t_m-\delta_m,t_m+\delta_m]\cap[0,T]} d_{i}(s).\label{def_c0im}
	\end{align}
	For $\epsArr\in\R^P_+$, we shall consider the following constraints
	\begin{align*}
	\V_0&:=\{\b x(\cdot)\in\Sx \,|\,\b C(t) \b x (t)\leq \b d(t),\,\forall \, t\in[0,T] \},\\
	\V_{\delta,\text{fin}}&:=\{\b x(\cdot)\in\Sx \,|\,\etaArr(\delta_m,t_m) \|\b x(\cdot)\|_K + \b C (t_m) \b x (t_m)\leq \b d(\delta_m,t_m),\, \forall \, m\in \iv{1}{N_0} \},\\
	\V_{\delta,\text{inf}}&:=\{\b x(\cdot)\in\Sx \,|\,\etaArr(\delta,t) \|\b x(\cdot)\|_K +\omeArr(\delta,t)+\b C(t) \b x (t)\leq \b d(t),\,\forall \, t\in[0,T]\},\\
	\V_{\epsilon}&:=\{\b x(\cdot)\in\Sx \,|\,\epsArr+\b C(t) \b x (t)\leq \b d(t),\,\forall \, t\in[0,T] \}.
	\end{align*}
	To these closed constraint sets correspond the problems $(\Psc_{0})$, $(\Psc_{\delta,\text{fin}})$ $(\Psc_{\delta,\text{inf}})$, and $(\Psc_{\epsilon})$. In particular,  $\bar{\b x}^0 (\cdot)$ denotes the optimal solution of $(\Psc_{0})$. \rv{When $\b C(\cdot)$ and $\b d(\cdot)$ are $\C^0$-continuous,} we prove right away that $\eta_i(\cdot, t)$ and $\omega_i(\cdot, t)$ converge uniformly in $t$ to $0$ as $\delta\rightarrow 0^+$, so that the SOC inequalities defining the set $\V_{\delta,\text{inf}}$ converge to the original affine constraints.
	\begin{Lemma}[Uniform continuity of $K$]\label{lem_modCo_noyau} If $\b A(\cdot)\in L^1(0,T)$ and $\b B(\cdot)\in L^2(0,T)$, then $K(\cdot,\cdot)$ is uniformly continuous. Assume furthermore that $\b C(\cdot)$ and $\b d(\cdot)$ are $\C^0$-continuous. Then for all $i\in\iv{1}{P}$, the increasing functions $\eta_i(\cdot,t)$ and $\omega_i(\cdot,t)$ converge to $0$ uniformly w.r.t.\ $t$ as $\delta\rightarrow 0^+$.
	\end{Lemma}		
	\noindent\textbf{Proof:} Since $\b A(\cdot)\in L^1(0,T)$, $\b \Phi_{\b A}(\cdot,\cdot)$ is uniformly continuous. Hence, for $\b Q\equiv\b 0$, through the explicit formulas \eqref{def_K0} and \eqref{def_K1}, we deduce that $K(\cdot,\cdot)$ is uniformly continuous. 
	Consequently, for any $i\in\iv{1}{P}$, and $s,t\in[0,T]$, recalling that $\BB_N$ denotes the closed unit ball of $\R^N$,
	\begin{align*}
	\|K(\cdot,t)\b c_{i}(t)-K(\cdot,s)\b c_{i}(s) \|_K\le  \|(K(\cdot,t)-K(\cdot,s))\b c_{i}(t) \|_K + \|K(\cdot,s)(\b c_{i}(t)-\b c_{i}(s)) \|_K&\\
	\le  \|\b c_{i}(t) \| \sup_{\b p\in\BB_N}\|(K(\cdot,t)-K(\cdot,s))\b p\|_{K} + \|\b c_{i}(t)-\b c_{i}(s) \|\sup_{\b p\in\BB_N}\|\b p^\top K(s,s)\b p\|^{1/2} &,
	\end{align*}
	which proves the statement for $\eta_i(\cdot,t)$, whereas the result for $\omega_i(\cdot,t)$ stems directly from the uniform continuity of $\b d(\cdot)$. Obviously, the components $\eta_i(\cdot,t)$ and $\omega_i(\cdot,t)$ are increasing for any given $t$.\\
	
	For $\b Q\not\equiv\b 0$, we do not \rv{have} explicit formulas such as \eqref{def_K0} and \eqref{def_K1}. Nonetheless, the $\|\cdot\|_K$-norm \eqref{def_1k-norm} for $\b Q\not\equiv\b 0$ is stronger than the $\|\cdot\|_K$-norm for $\b Q\equiv\b 0$. Since, for $\b Q\equiv\b 0$, $K$ is uniformly continuous, owing to \cite[Proposition 24]{Schwartz1964}, the topology induced by $K$ over $\Sx$ is stronger than the topology of uniform convergence over $[0,T]$. Hence the topology induced by $K$ for $\b Q\not\equiv\b 0$ is also stronger. Therefore, using again the result of \citet{Schwartz1964}, for $\b Q\not\equiv\b 0$, $K(\cdot,\cdot)$ is continuous w.r.t.\ each variable  and locally bounded.\footnote{This allows to derive the continuity of $\eta_i(\cdot,t)$ but does not provide a uniform bound w.r.t to $t$.} Hence $K(\cdot,\cdot)$ is bounded on the compact set $[0,T]\times[0,T]$. Let us prove the continuity of $t\mapsto K(t,t)$. Since $K$ is bounded, by \eqref{eq_integralK}, $\tilde{K}$ is bounded, so $\b U_t(\cdot)\in L^2(0,T)$. Let $\b p\in\BB_N$ and $t\in[0,T]$, then, by definition of $K$, $K(\cdot,t)\b p \in \Sx$ is associated to the control $\b U_t(\cdot)\b p$. Let $\delta>0$, by \rv{the variation of constants \rvS{formula}} \eqref{eq_Duhamel}, we have
	\begin{align*}
	\|(K(t+\delta,t)- K(t,t))\b p\|&\le  \|\b \Phi_{\b A}(t+\delta,0)-\Phi_{\b A}(t,0) \| \cdot \|K(0,t)\b p\|
	+\int_{t}^{t+\delta} \|\b \Phi_{\b A}(t,s) \b B(s) \b U_t(s)\b p \| \d s.
	\end{align*}
	Let $\lambda>0$. With a similar computation when permuting $t$ and $t+\delta$, taking the supremum over $\b p\in\BB_N$, one can find $\Delta>0$ such that for any $\delta\in[0,\Delta]$,
	$$ \max(\|K(t+\delta,t)- K(t,t)\|, \|K(t,t+\delta)- K(t+\delta,t+\delta)\|)\le \lambda/2.$$
	Hence, owing to the Hermitian symmetry of $K$, for any $\delta\in[0,\Delta]$,
	$$\|K(t+\delta,t+\delta)- K(t,t)\|\le \|K(t+\delta,t)- K(t,t)\|+\|K(t,t+\delta)- K(t+\delta,t+\delta)\|\le \lambda.$$
	This shows that $K$ is indeed continuous on the diagonal.	As underlined by Laurent Schwartz, showing the continuity of $t\mapsto K(t,t)$ is enough to conclude. We reproduce briefly his argument \citep[see][p194]{Schwartz1964}: for any $\b p\in\R^N$, whenever $t$ converges to $t_0$, $K(\cdot,t)\b p$ weakly converges in $\Sx$ to $K(\cdot,t_0)\b p$, however, by continuity on the diagonal, the norm $\b p^\top K(t,t)\b p$ converges to $\b p^\top K(t_0,t_0)\b p$, so $K(\cdot,t)\b p$ strongly converges in $\Sx$ and by extension in $\C^0(0,T)$, exactly showing that $K(\cdot,\cdot)$ is continuous, hence uniformly continuous.
	\begin{flushright}
		$\blacksquare$
	\end{flushright}

	Generically, under the minimal assumptions of Lemma \ref{lem_modCo_noyau}, one can apply the SOC-scheme to obtain $(\Psc_{\delta,\text{fin}})$, which is equivalent to a finite dimensional problem owing to the representer theorem (Theorem \ref{thm_representer}). By Lemma \ref{lem_modCo_noyau}, the scheme is coherent since for $\delta$ decreasing to zero, the coefficients of the SOC constraints converge uniformly in $t$ to those of the original problem $(\Psc_{0})$. However, 
	ensuring that the solution $\bar{\b x}^{\delta,\text{fin}}(\cdot)$ converges to $\bar{\b x}^{0}(\cdot)$ requires a more thorough analysis.
	
	\begin{Proposition}[Nested sequence]\label{thm_incl_SOC} Let $\delta_{\text{max}}:=\max_{m\in \iv{1}{N_0} }\delta_m$. For any $\delta \geq \delta_{\text{max}}$, if, for a given $y_0\geq 0$, $\epsArr \geq \sup_{t\in[0,T]} [\etaArr(\delta,t) y_0 +\omeArr(\delta,t)]$, then we have a nested sequence 
		\begin{equation}\label{eq_nested_sets}
		(\V_{\epsilon}\cap y_0\BB_K)\subset \V_{\delta,\text{inf}}\subset \V_{\delta,\text{fin}} \subset \V_0.
		\end{equation}
	\end{Proposition}		
	\noindent\textbf{Proof:} The inclusion $\V_{\epsilon}\cap y_0\BB_K\subset \V_{\delta,\text{inf}}$ stems from the definition of the sets. Since $d_{i}(\delta_m,t_m) \geq d_{i}(t_m) - \omega_i(\delta_m,t_m)$,  $\V_{\delta,\text{inf}}\subset \V_{\delta,\text{fin}}$. Recall that $[0,T]\subset\cup_{m\in\iv{1}{N_0}} [t_m-\delta_m,t_m+\delta_m]$. Let $t\in[0,T]$ and $\b x(\cdot)\in\V_{\delta,\text{fin}}$. Take $m\in\iv{1}{N_0}$ such that $t\in [t_m-\delta_m,t_m+\delta_m]$. For any $i\in\iv{1}{P}$, applying the reproducing property and Cauchy-Schwarz inequality, 
	\begin{align*}
		\b c_i(t)^\top \b x (t) &=  \b c_i (t_m)^\top \b x (t_m) + \langle \b x(\cdot), K(\cdot,t)\b c_{i}(t)-K(\cdot,t_m)\b c_{i}(t_m) \rangle_K \\
		\b c_i(t)^\top \b x (t) &\le \b c_i (t_m)^\top \b x (t_m) + \eta_i(\delta_m,t_m) \|\b x(\cdot)\|_K \le d_i(\delta_m,t_m) \le d_i(t),
	\end{align*}
	$d_i(\delta_m,t_m)$ being by definition the infimum of the $i$-th component $d_i(t)$ of $\b d(t)$ on $[t_m-\delta_m,t_m+\delta_m]$. So $\b x(\cdot)\in \V_0$, hence $\V_{\delta,\text{fin}} \subset \V_0$.
		\begin{flushright}
			$\blacksquare$
		\end{flushright}
	Proposition \ref{thm_incl_SOC} states that, on the one hand, enforcing a finite number of SOC constraints with $\etaArr$ as in \eqref{def_eta} is more restrictive than enforcing an infinite number of affine constraints. On the other hand, SOC constraints are less restrictive than shrinking the affine constraints by some $\epsArr>0$. The nested property \eqref{eq_nested_sets} is instrumental in our analysis. As a matter of fact, we shall focus on $\epsilon$-perturbations of affine constraints rather than on SOC constraints to construct a trajectory $\b x^\epsilon(\cdot)\in\V_{\epsilon}$ close to $\bar{\b x}^0(\cdot)$. We shall then resort to strong convexity arguments to derive bounds on $\|\bar{\b x}^{\delta,\text{fin}}(\cdot)-\bar{\b x}^{0}(\cdot)\|_K$.\\

	We now list the hypotheses used to prove our main result. 
	\begin{description}
		\item  [(H-gen)\label{hyp_general}] $\b A(\cdot)\in L^1(0,T)$ and $\b B(\cdot)\in L^2(0,T)$, $\b C(\cdot)$ and $\b d(\cdot)$ are $\C^0$-continuous. 
		\item  [(H-sol)\label{hyp_existSol}] $\b C(0)\b x_0< \b d(0)$ and there exists $\epsArr>0$ such that $\V_\epsilon\cap\{\b x(\cdot)\,|\,\b x(0)= \b x_0\}\neq\emptyset$, i.e.\ there exists a trajectory $\b x^{\epsilon}(\cdot)\in\Sx$ satisfying strictly the affine constraints, as well as the initial condition, with $\b x_0$ interior to the state constraints.
		\item  [(H-L)\label{hyp_str_convexity}] There exists $\mu>0$ such that the objective function $\Lcal:\b x(\cdot)\in\Sx \mapsto g(\b x(T))+\|\b x(\cdot)\|^2_K$ is $\mu$-strongly convex. The terminal cost $g(\cdot)$ is continuous over $\R^N$, $\b Q(\cdot)\in L^1(0,T)$, and $\b R(\cdot)\in L^2(0,T)$. There exists $r>0$ such that $\b R(t) \succcurlyeq r \Id_M$ for all $t\in[0,T]$.
	\end{description} 

	\noindent\textbf{Discussion of the Assumptions:} Assumption \ref{hyp_general} ensures the $\C^0$-continuity of the kernel $K$ and of the functions $\eta_i$ and $\omega_i$ (Lemma \ref{lem_modCo_noyau}). Assumption \ref{hyp_str_convexity} concerns the objective function, whereas \ref{hyp_existSol} ensures that the set of trajectories satisfying the state constraints is non-empty if the latter are shrunk:
	
	\begin{itemize}[labelindent=0cm,leftmargin=*,topsep=0cm,partopsep=0cm,parsep=0.1cm,itemsep=0cm]
		\item The existence requirement in \ref{hyp_existSol} can be derived from assumptions on the existence of interior viable trajectories. We provide in the Annex an example of such assumptions (Lemma \ref{lem_IPC_cndViab}) based on inward pointing conditions on the boundary and on regularity assumptions on the constraints and the dynamics. \rvS{The assumption $\b C (0) \b x_0 < \b d_0$ ensures that $\b x_0$ is a suitable initial condition for the $\epsArr$-tightening $\V_{\epsilon}$.}\footnote{\rvS{Since the SOC tightening lies in-between the $\epsArr$-tightening and the original constraints (Proposition \ref{thm_incl_SOC}), it cannot be guaranteed that an initial condition on the border of the constraints would be suitable for the SOC tightening.}}  
		\item The strong convexity requirement in \ref{hyp_str_convexity} is obviously satisfied whenever $g(\cdot)$ is convex.\footnote{More generally, $g$ could be $\mu_0$-semiconvex (i.e.\ $g(\cdot)+\frac{\mu_0}{2}\|\cdot\|^2$ is convex) with $2>\mu_0\sup_{\b p\in\BB_N}\|\b p^\top K(T,T)\b p\|^{1/2}$.} It is required in order to bound the distance on solutions since the problems $(\Psc_{*})$ share the same objective but different constraint sets. The terminal cost $g(\cdot)$ is supposed $\C^0$-continuous, but could be taken merely locally continuous in a neighborhood\footnote{This neighborhood is considered with respect to the relative topology of the terminal constraint set $\{\b x\in\R^N \,|\,\b C(T) \b x\leq \b d(T)\}.$} of $\bar{\b x}^0(T)$ and lower bounded over any compact subset of $\R^N$. By Lemma \ref{lem_existenceSol}, $(\Psc_{0})$ has a unique optimal solution $\bar{\b x}^0(\cdot)$.
	\end{itemize}

	\begin{Lemma}[Existence and uniqueness of solutions]\label{lem_existenceSol} Under Assumptions \ref{hyp_general}, \ref{hyp_existSol}, and \ref{hyp_str_convexity}, $\bar{\b x}^0(\cdot)$ exists and is unique.	The same result holds true for $(\Psc_{\delta,\text{fin}})$, $(\Psc_{\delta,\text{inf}})$, and $(\Psc_{\epsilon})$ for any $\delta\in[0,\delta_0]$, where $\delta_0>0$ satisfies that $\epsArr\ge \sup_{t\in[0,T]}  [\etaArr(\delta_0,t) \|\b x^{\epsilon}(\cdot)\|_K +\omeArr(\delta_0,t)]$ for $\epsArr$ as in \ref{hyp_existSol}.
	\end{Lemma}	
	\noindent\textbf{Proof:} The existence result is a consequence of Tonelli's direct method, usually stated for lower bounded and lower semi-continuous $g(\cdot)$. We detail the proof since our Assumptions are both slightly different and stronger. Since by Assumption \ref{hyp_existSol}, $\V_0\neq\emptyset$, let $(\b x_n(\cdot),\b u_n(\cdot))$ be a minimizing sequence of \eqref{opt-LQR} converging to the optimal value $\bar{m}$. As $\Lcal(\cdot)$ is $\mu$-strongly convex, $\frac{\Lcal(\b x(\cdot))}{\|\b x(\cdot)\|_K}\rightarrow +\infty$ as $\|\b x(\cdot)\|_K\rightarrow +\infty$. Hence $\bar{m}$ is finite, and $(\b x_n(\cdot))_n$ is a subset of a ball $M_0\BB_K\subset\Sx$, for some $M_0>0$. Since $\b x(\cdot)\in \Sx\mapsto \b x(T)$ is continuous, $\{\b x(T)\,|\, \b x(\cdot) \in M_0\BB_K\}$ is also bounded. By continuity of $g(\cdot)$, let $m_g:=-\inf_{\b x(\cdot) \in M_0\BB_K}g(\b x(T))<+\infty$. 
	Consequently, for $n$ large enough,
	\begin{align*}
	r \Vert\b u_n(\cdot)\Vert^2_{L^2(0,T)} \le \Vert\b x_n(\cdot)\Vert^2_K \le \bar{m} + 1 - g(\b x_n(T)) \le \bar{m} + 1 +m_g,
	\end{align*}
	so $(\b u_n(\cdot))_n$ is bounded in $L^2$, and we can take a subsequence $(\b u_{n_i}(\cdot))_i$  weakly converging to some $\b u(\cdot)$. Let $s,t\in[0,T]$. By \rv{the variation of constants \rvS{formula}} \eqref{eq_Duhamel}, we had derived \eqref{eq_Duhamel_diff}.	Since $\b A(\cdot)\in L^1(0,T)$ and $\b B(\cdot)\in L^2(0,T)$, and $\b x_n(s)$ is uniformly bounded in $n$ and $s$ by the reproducing property and continuity of $K$ (Lemma \ref{lem_modCo_noyau}), $(\b x_{n_i}(\cdot))_i$ is equicontinuous. By Ascoli's theorem, we thus have a subsequence $(\b x_{n_{i_j}}(\cdot))_j$ uniformly converging to some $\b x(\cdot)$ satisfying \eqref{eq_Duhamel} for $\b u(\cdot)$, thus $\b x(\cdot)\in\Sx$. By continuity of $g(\cdot)$, $\Lcal(\b x(\cdot))=\bar{m}$. Since $\Lcal(\cdot)$ is strongly convex, the optimal trajectory is unique and belongs to the closed set $\V_0$. To conclude, replace $\V_0$ with $\V_{\delta,\text{fin}}$ (resp.\ $\V_{\delta,\text{inf}}$, and $\V_{\epsilon}$), the inequality satisfied by $\delta_0$ shows that  $\b x^{\epsilon}(\cdot)\in\V_{\delta,\text{fin}}$, \rv{consequently the constraint sets} are non-empty. The same arguments as above yield the result. 
	\begin{flushright}
		$\blacksquare$
	\end{flushright}
	\begin{Theorem}[Main result - Approximation by SOC constraints]\label{thm_bounds_SOC}
		Under Assumptions \ref{hyp_general}, \ref{hyp_existSol}, and \ref{hyp_str_convexity}, for any $\lambda>0$, there exists $\bar{\delta}>0$ such that for all $N_0>0$ and $(\delta_m)_{m\in\iv{1}{N_0}}$, with $[0,T]\subset\cup_{m\in\iv{1}{N_0}} [t_m-\delta_m,t_m+\delta_m]$ satisfying $\bar{\delta}\ge\max_{m\in \iv{1}{N_0} }\delta_m$, we have
		\begin{equation}\label{eq_bounds_SOC}
		\frac{1}{\gamma_K}\cdot\sup_{t\in[0,T]}\|\bar{\b x}^{\delta,\text{fin}}(t)-\bar{\b x}^{0}(t)\|	\le \|\bar{\b x}^{\delta,\text{fin}}(\cdot)-\bar{\b x}^{0}(\cdot)\|_K \le \lambda
		\end{equation}
		with $\gamma_K:=\sup_{t\in[0,T], \,\b p\in\BB_N}\sqrt{\b p^\top K(t,t) \b p}$.
	\end{Theorem}
	\noindent\textbf{Proof:} Let $\lambda>0$. Consider any $\tilde{\lambda}>0$ such that 
	\begin{equation}\label{eq_tech_lambda}
		2\tilde{\lambda}+\tilde{\lambda}(\tilde{\lambda}+2\|\bar{\b x}^0(\cdot)\|_{L^{\infty}\left(0,T\right)})\|\b Q(\cdot)\|_{L^{1}(0,T)}\le \lambda.
	\end{equation}
	By Assumption \ref{hyp_existSol}, pick $\b x^\epsilon(\cdot)\in\V_\epsilon$ such that $\b x^\epsilon(0)=\bar{\b x}^0(0)$. Denote by $\b u^\epsilon(\cdot)$ the associated control. Take $\alpha>0$ small enough such that $\b x^{\alpha\epsilon}(\cdot):= \alpha \b x^\epsilon(\cdot)+ (1-\alpha)\bar{\b x}^0(\cdot)\in\Sx$ and $\b u^{\alpha\epsilon}(\cdot):= \alpha \b u^\epsilon(\cdot)+ (1-\alpha)\bar{\b u}^0(\cdot)$ satisfy
	\begin{align*}
	\|\bar{\b x}^0(\cdot)-\b x^{\alpha\epsilon}(\cdot)\|_{L^{\infty}\left(0,T\right)}=\alpha \|\bar{\b x}^0(\cdot)-\b x^{\epsilon}(\cdot)\|_{L^{\infty}\left(0,T\right)} &\leq \tilde{\lambda}\\
	\vert\Vert\b R(\cdot)^{1/2}\bar{\b u}(\cdot)\Vert^2_{L^2(0,T)}-\Vert \b R(\cdot)^{1/2}\b u^{\alpha\epsilon}(\cdot)\Vert^2_{L^2(0,T)}\vert&\leq\tilde{\lambda}
	\end{align*}
	and, by continuity of $g(\cdot)$, $|g(\bar{\b x}^0(T))-g(\b x^{\alpha\epsilon}(T))|\leq\tilde{\lambda}$. Consequently $\b x^{\alpha\epsilon}(0)=\bar{\b x}^0(0)$ and for all $t\in[0,T]$, $\b C(t) \b x^{\alpha\epsilon}(t) \le \alpha (\b d(t)-\epsArr) + (1-\alpha)\b d(t)= \b d(t)-\alpha\epsArr$,
	so $\b x^{\alpha\epsilon}(\cdot)\in\V_{\alpha\epsilon}$. Hence
	\begin{align*}
		\Lcal(\b x^{\alpha\epsilon}(\cdot))-\Lcal(\bar{\b x}^0(\cdot))&\le |g(\bar{\b x}^0(T))-g(\b x^{\alpha\epsilon}(T))|+\left\vert\|\bar{\b x}^0(\cdot)\|^2_K -\|\b x^{\alpha\epsilon}(\cdot)\|^2_K\right\vert\\
		&\hspace{-3.5cm}\le \tilde{\lambda} + \int_{0}^{T}\left\vert(\bar{\b x}^0(t)-\b x^{\alpha\epsilon}(t))^{\top}\b Q(t) (\bar{\b x}^0(t)+\b x^{\alpha\epsilon}(t))\right\vert \d t + \left\vert\Vert\b R(\cdot)^{1/2}\bar{\b u}(\cdot)\Vert^2_{L^2(0,T)}-\Vert\b R(\cdot)^{1/2}\b u^{\alpha\epsilon}(\cdot)\Vert^2_{L^2(0,T)}\right\vert\\
		&\hspace{-3.5cm}\le 2 \tilde{\lambda} + \tilde{\lambda}(\tilde{\lambda}+2\|\bar{\b x}^0(\cdot)\|_{L^{\infty}\left(0,T\right)})\|\b Q(\cdot)\|_{L^{1}(0,T)}\stackrel{\eqref{eq_tech_lambda}}{\le} \lambda.
	\end{align*}
	Let $\delta_0>0$ such that $\alpha\epsArr\ge \sup_{t\in[0,T]}  [\etaArr(\delta_0,t) \|\b x^{\alpha\epsilon}(\cdot)\|_K +\omeArr(\delta_0,t)]$. Then $\b x^{\alpha\epsilon}(\cdot)\in \V_{\delta_0,\text{inf}}\subset\V_{\delta_0,\text{fin}}$, the sets thus being non-empty. Notice that, for any $\delta\in[0,\delta_0]$, as $\bar{\b x}^{*}(\cdot)$ is optimal for $(\Psc_{*})$, from the nested property \eqref{eq_nested_sets}, we derive that
	\begin{equation*}
		\Lcal( \bar{\b x}^{\delta,\text{fin}}(\cdot))\le \Lcal( \bar{\b x}^{\delta,\text{inf}}(\cdot))\le \Lcal(\bar{\b x}^{\delta_0,\text{inf}}(\cdot)).
	\end{equation*}
	As $\Lcal(\cdot)$ is $\mu$-convex, $\Lcal^{-1}(]-\infty, \Lcal( \bar{\b x}^{\delta_0,\text{inf}}(\cdot))])$ is a bounded set, contained in a ball $M_0\BB_K$ for some $M_0>0$, and containing all the $\{\bar{\b x}^{\delta,\text{inf}}(\cdot)\}_{\delta\in[0,\delta_0]}$. Since $\Sx$ is a vRKHS, $\b x(\cdot)\in \Sx\mapsto \b x(T)$ is continuous. So $\{\b x(T)\,|\, \b x(\cdot) \in M_0\BB_K\}$ is also bounded. Hence, for any $\delta\in[0,\delta_0]$,
	\begin{gather*}
	g(\bar{\b x}^{\delta,\text{inf}}(T)) + \|\bar{\b x}^{\delta,\text{inf}}(\cdot)\|^2_K \le \Lcal(\bar{\b x}^{\delta_0,\text{inf}}(\cdot)) \le |g(\bar{\b x}^{\delta_0,\text{inf}}(T))| +\|\bar{\b x}^{\delta_0,\text{inf}}(\cdot)\|^2_K\\
		\|\bar{\b x}^{\delta,\text{inf}}(\cdot)\|_K \le \|\bar{\b x}^{\delta_0,\text{inf}}(\cdot)\|_K +\sqrt{|g(\bar{\b x}^{\delta_0,\text{inf}}(T))|} + |\inf_{\b x(\cdot) \in M_0\BB_K}g(\b x(T))|^{\frac{1}{2}}=:y_0.
	\end{gather*}
	As $\|\bar{\b x}^{\delta_0,\text{inf}}(\cdot)\|_K\le y_0$, $\b x^{\alpha\epsilon}(\cdot)$ and $\bar{\b x}^{\delta_0,\text{inf}}(\cdot)$ are both admissible for the following problem, 
	$$\min_{\substack{\b x(\cdot)\in\V_{\delta_0,\text{inf}}\\
				\|\b x(\cdot)\|_K\le y_0+\|\bar{\b x}^{\alpha\epsilon}(\cdot)\|_K}}\Lcal(\b x(\cdot))$$
	with $\bar{\b x}^{\delta_0,\text{inf}}(\cdot)$ being optimal by definition, hence we have $\Lcal(\bar{\b x}^{\delta_0,\text{inf}}(\cdot))\le \Lcal(\b x^{\alpha\epsilon}(\cdot))$. To conclude, let $\bar{\delta}\in]0,\delta_0]$ such that $\alpha\epsArr\ge \sup_{t\in[0,T]}  [\etaArr(\bar{\delta},t) y_0 +\omeArr(\bar{\delta},t)]$. Then, for any $\delta\in]0,\bar{\delta}]$, by strong convexity of $\Lcal(\cdot)$, $\bar{\b x}^{0}(\cdot)$ being optimal for $(\Psc_0)$,
	\begin{align*}
		\frac{\mu}{2} \|\bar{\b x}^{\delta,\text{fin}}(\cdot)-\bar{\b x}^{0}(\cdot)\|_K^2 \le \Lcal( \bar{\b x}^{\delta,\text{fin}}(\cdot))-\Lcal( \bar{\b x}^{0}(\cdot))\le \Lcal(\b x^{\alpha\epsilon}(\cdot))-\Lcal( \bar{\b x}^{0})\le \lambda.
	\end{align*}
	Replacing $\lambda$ by $\sqrt{2\lambda/\mu}$, we deduced that $\|\bar{\b x}^{\delta,\text{fin}}(\cdot)-\bar{\b x}^{0}(\cdot)\|_K \le \lambda$. By Cauchy-Schwarz inequality, for any $t\in[0,T]$, $\|\bar{\b x}^{\delta,\text{fin}}(t)-\bar{\b x}^{0}(t)\|	\le \|\bar{\b x}^{\delta,\text{fin}}(\cdot)-\bar{\b x}^{0}(\cdot)\|_K \sup_{\,\b p\in\BB_N}\sqrt{\b p^\top K(t,t) \b p}.$
	By definition of $\gamma_K$, taking the supremum over $[0,T]$, we derive the remaining inequality.
	\begin{flushright}
		$\blacksquare$
	\end{flushright}
	Theorem \ref{thm_bounds_SOC} states that, when the discretization steps $(\delta_{\text{m}})_m$ go to zero, then the solution $\bar{\b x}^{\delta,\text{fin}}(\cdot)$ of the SOC-approximation $(\Psc_{\delta,\text{fin}})$ can be made arbitrarily close to the solution $\bar{\b x}^0(\cdot)$ of the original problem $(\Psc_{0})$, uniqueness being ensured by Assumption \ref{hyp_str_convexity}. Concerning the stability of solutions under shrinking perturbation of the state constraints, the result of Theorem \ref{thm_bounds_SOC} actually also holds when replacing $\bar{\b x}^{\delta,\text{fin}}(\cdot)$ by $\bar{\b x}^{\epsilon}(\cdot)$, showing that, when $\epsArr$ goes to zero, $\bar{\b x}^{\epsilon}(\cdot)$ converges to $\bar{\b x}^0(\cdot)$. 

	\section{Finite-dimensional implementation and numerical example}\label{sec_numerics}
	In this section, we express the finite-dimensional equivalent of problem \eqref{opt-LQR_SOC} owing to the representer theorem (Theorem \ref{thm_representer}) and discuss its implementation on a numerical example. In general, the SOC transformation requires only the minimal hypotheses of Lemma \ref{lem_modCo_noyau} to be conceptually grounded. Theorem \ref{thm_bounds_SOC} \rv{essentially states theoretical guarantees of convergence for small discretization steps}. For numerical applications, the problem can thus be extended to incorporate costs or equality constraints at any finite number of intermediate times, as in path planning problems. This situation was already met when discussing the controllability Gramian \eqref{opt-LQ_Gram}. We consequently enrich \eqref{opt-LQR_SOC} to optimization problems considered in Theorem \ref{thm_representer}, of the following form, with $\|\b x(\cdot)\|^2_{K}=\|\b x(0)\|^2+\|\b x(\cdot)\|^2_{K_1}$, 
	\begin{mini}
		{\substack{\b x(\cdot)\in\Sx}}{L\left(\b c_{0,1}^{\top}\b x(t_{0,1}),\dots, \b c_{0,N_0}^{\top}\b x(t_{0,N_0})\right)+\|\b x(\cdot)\|^2_{K}\phantom{\eta_i(\delta_{i,m},t_{i,m}) \|\b x(\cdot)\|_K + \b c_i (t_{i,m})^{\top} \b x (t_{i,m})\leq d_{i} }}{\label{opt-LQR_SOC_pathPlan}\tag{$\Psc_{\text{SOC}}$}}{}
		\addConstraint{ \hspace{-5mm}\eta_i(\delta_{i,m},t_{i,m}) \|\b x(\cdot)\|_K + \b c_i (t_{i,m})^{\top} \b x (t_{i,m})}{\leq d_{i}(\delta_{i,m},t_{i,m}),\, \forall \, m\in \iv{1}{N_i}, \forall \,i\in\iv{1}{P},}
	\end{mini}
	for $\{\b c_{0,m}\}_{m\in\iv{1}{N_0}} \subset\R^N$ . The differences between \eqref{opt-LQR_SOC} and \eqref{opt-LQR_SOC_pathPlan} are that $g(\cdot)$, which \rv{depended only} on the terminal point, is now replaced with a loss $L:\RR^{N_0}\rightarrow \RR\cup\{+\infty\}$, defined on a finite number $(t_{0,m})_{m\in\iv{1}{N_0}}$ of intermediate points (taken with repetition), and that different discretization grids $(\delta_{i,m},t_{i,m})_{i,m\in\iv{1}{N_i}}$ are used for each constraint $i\in\iv{1}{P}$. As no structural assumptions are imposed on $L$, it may incorporate indicator functions to account for the initial condition or for rendezvous points.\footnote{To write \eqref{opt-LQR_SOC} as \eqref{opt-LQR_SOC_pathPlan}, take $N_0=2N$, $t_{0,1}=\dots=t_{0,N}=0$,  $t_{0,N+1}=\dots=t_{0,N_0}=T$. Denoting by $\chi_{\b x_0}(\cdot)$ the indicator function of $\b x_0$, set $L\left(\b c_{0,1}^{\top}\b x(t_{0,1}),\dots, \b c_{0,N_0}^{\top}\b x(t_{0,N_0})\right):= \chi_{\b x_0}(\b x(t_{0,1}))+ g(\b x(t_{0,N_0}))$.} In Section \ref{sec_bounds}, the SOC constraints \rv{were introduced} to turn an infinite number of affine constraints over $[0,T]$ into a finite number of SOC constraints. The logic is therefore to separate the discrete pointwise requirements (which go to $L$) from the constraints that should hold on $[0,T]$ (which are approximated by SOC constraints). Since the constraints on $[0,T]$ may apply to different components of the state, we may consider different grid steps for each $i$.
	
	In general, adding the SOC terms leads to more conservative solutions w.r.t.\ the one with affine constraints. \rvS{The SOC constraints formulation only requires estimating $\eta_i(\delta_{i,m},t_{i,m})$ and $-d_{i}(\delta_{i,m},t_{i,m})$ defined in \eqref{def_eta} and \eqref{def_c0im}. Both quantities should in principle be overestimated for the guarantees to hold. However, since the tightening results from a worst-case Cauchy-Schwarz scenario as discussed in \eqref{eq:mod_continuity_constraint}, in practice, underestimation does not affect the numerical results.} Moreover the definition \eqref{def_eta} of $\eta_i$ is not the only possible formulation.\footnote{Using that $\b x(\cdot)=\b z_0(\cdot)+\b z_1(\cdot)\in\Sx_0\oplus\Sx_u$ one could consider two $\eta$-terms instead of one to derive a less conservative tightening. In the same spirit $\b d(t)$ could be projected onto $\Sx$, and the projection incorporated in the scalar products of the left hand-side. For $\b Q(\cdot)\not\equiv \b 0$, $\|\b x(0)\|^2$ could be replaced in \eqref{def_1k-norm} by an $\Sx_0$-norm $\|\b S_0 \b x(0)\|^2$ with a surjective $\b S_0 \in \R^{N_0,N}$ where $N_0 = \text{dim}(\Sx_0)$. This would not change the formulation, but lead to a ``tighter`` norm $\|\b x(\cdot)\|_K$.} Its choice results from geometrical considerations on coverings\footnote{Definition \eqref{def_eta} corresponds to a covering made of balls in $\Sx$.} of compact sets in infinite-dimensional Hilbert spaces (see Section 3 of \citet{aubin2020hard_SDP} for more details). Even for other values of $\eta_i$ than \eqref{def_eta}, considering SOC terms in the constraints proves to be beneficial in terms of local satisfaction of the constraints on a neighborhood of $t_{i,m}$. Besides, the discretization grids considered here are 'static' in the sense that they are fixed before solving \eqref{opt-LQR_SOC_pathPlan}. Extensions to 'dynamic' grids, refined depending on the optimization steps, can be found in Section 5 of \citet{aubin2020hard_SDP}. 
	
	By Theorem \ref{thm_representer}, for $\b x(\cdot)= \sum_{j=0}^{P} \sum_{m=1}^{N_j} K(\cdot,t_{j,m}) \b p_{j,m}$ and $z=\|\b x(\cdot)\|_{K}$,\footnote{\rv{The reproducing property, $\langle K(\cdot,t_1)\b p_1,K(\cdot,t_2)\b p_2\rangle_K=\b p_2^\top K(t_2,t_1)\b p_1$, applied to $\b x(\cdot)$ allows to explicit $z$.}} \eqref{opt-LQR_SOC_pathPlan} is equivalent to 
	\begin{mini*}
		{\substack{z\in\R_+,\\j\in\iv{0}{P},\,m\in  \iv{1}{N_j},\\ \b p_{j,m}\in\R^{N},\; \alpha_{j,m} \in\R}}{L\left(\left(\sum_{j=0}^{P} \sum_{m=1}^{N_j} K(t_{0,n},t_{j,m}) \b p_{j,m}\right)_{n\in\iv{1}{N_0}}\right)+ z^2}{}{}
		\addConstraint{\hspace{-8cm}z^2=\sum_{i=0}^{P} \sum_{n=1}^{N_i}\sum_{j=0}^{P} \sum_{m=1}^{N_j} \b p_{i,n}^\top K(t_{i,n},t_{j,m}) \b p_{j,m}}{}{}
		\addConstraint{\hspace{-8cm} \b p_{j,m}=\alpha_{j,m} \b c_j(t_m),\quad \forall \, m\in \iv{1}{N_j}, \forall \,j\in\iv{1}{P}}{}
		\addConstraint{\let\scriptstyle\textstyle\substack{ \hspace{-1cm}\eta_i(\delta_{i,m},t_{i,m})z + \sum_{j=0}^{P}\sum_{m=1}^{N_j} \b c_i (t_{i,m})^{\top} K(t_{i,m},t_{j,m}) \b p_{j,m}\\  \hspace{-1cm}\;\leq d_{i}(\delta_{i,m},t_{i,m}),}}{\quad \forall \, m\in \iv{1}{N_i}, \forall \,i\in\iv{1}{P}}{}.
	\end{mini*}
	\rvS{For quadratic $L$ with indicator functions, after incorporating the quadratic objective into the constraints through a change of variables, \eqref{opt-LQR_SOC_pathPlan} writes as a second-order cone program (SOCP).} Hence this problem can be straightforwardly implemented in convex solvers or modeling frameworks such as CVXGEN \citep{mattingley2012CVX}. SOCP is slightly more expensive computationally than the quadratic programs (QP) classically derived for LQR \citep{kojima2004lq}.  \rvS{Both QP and SOCP have polynomial computation times, however the exponent is problem-dependent, so achieving a theoretical comparison is challenging. From a numerical viewpoint, it seems to be also very much framework-dependent. For instance, when using CVXGEN, 
	the baseline cost of calling CVX blurs the difference in computation times between QP and SOCP.} As the more time points, the more coefficients, it is beneficial to define the grids as subsets of a 'master grid'. \rv{Furthermore when computing $\eta_i$ or $K(s,t)$, one has to approximate the supremum in \eqref{def_eta} or the integral in \eqref{def_K1}, for instance through sampling. For time-invariant dynamics and $\b Q(\cdot) \equiv \b 0$, one can use the method proposed in \citet{VanLoan1978computingExpM} to quickly compute $K_1(s,t)$ in \eqref{def_K1}.} \rvS{If sought for, the analysis of the approximation error would correspond to that of the stability of the solution when perturbing the constants and matrices appearing in the SOCP. This leads to a large body of articles (e.g.\ \citet{Bonnans2005} and quoting articles) which exceeds the purposes of this study. If one overestimates the value of $\eta$ and adds a term $\lambda_{\text{cond}}\sum_{i=0}^{P} \sum_{n=1}^{N_i} \| \b p_{i,n}\|^2$ to $z^2$, then the state constraints are further tightened. This is done in practice to improve the conditioning of the matrices, and discussed along other numerical implementation details of a SOC-constrained kernel regression in \citet{aubin2020hard_nips}. For a sufficiently small error in $K$, constraints satisfaction could be guaranteed, even when facing numerical errors, through this further tightening involving $\lambda_{\text{cond}}$.}\\ 

	To highlight the behavior of the SOC-transformation \eqref{opt-LQR_SOC_pathPlan} of problem \eqref{opt-LQR}, we consider the problem \eqref{opt-pendulum} of a 'linear' pendulum with angle $x(t)$ where we control the derivative $u(t)$ of a \rv{forcing term} $w(t)$, with state constraints both on $w(t)$ and on $\dot{ x}(t)$, \rv{the full state being $\b x:=[ x, \dot{ x}, w]\in\R^3$,}
	\begin{mini}
	{\substack{ x(\cdot), u(\cdot)}}{-\lambda_T\, \dot{ x}(T)+ \lambda_u \Vert u(\cdot)\Vert_{L^2(0,T)}^2 + \|\b x(0)\|^2 }{\label{opt-pendulum}\tag{$\Psc_{\text{pend}}$}}{}
	\addConstraint{  x(0)=0.5, \quad \dot{ x}(0)=0,\quad  w(0)=0 }{}
	\addConstraint{x(T/3)=0.5,\quad   x(T)=0}{}
	\addConstraint{ \ddot{ x}(t)=-10 \,  x(t) + w(t), \quad  \dot{w}(t)=u(t), \, \text{a.e.\ in} \, [0,T]}{}
	\addConstraint{ \dot{ x}(t)\in [-3,+\infty[, \quad w(t)\in[-10,10],\,\forall \, t\in[0,T]}{}.
	\end{mini}
\begin{figure}
	\includegraphics[keepaspectratio, width=\textwidth,height=.75\textheight]{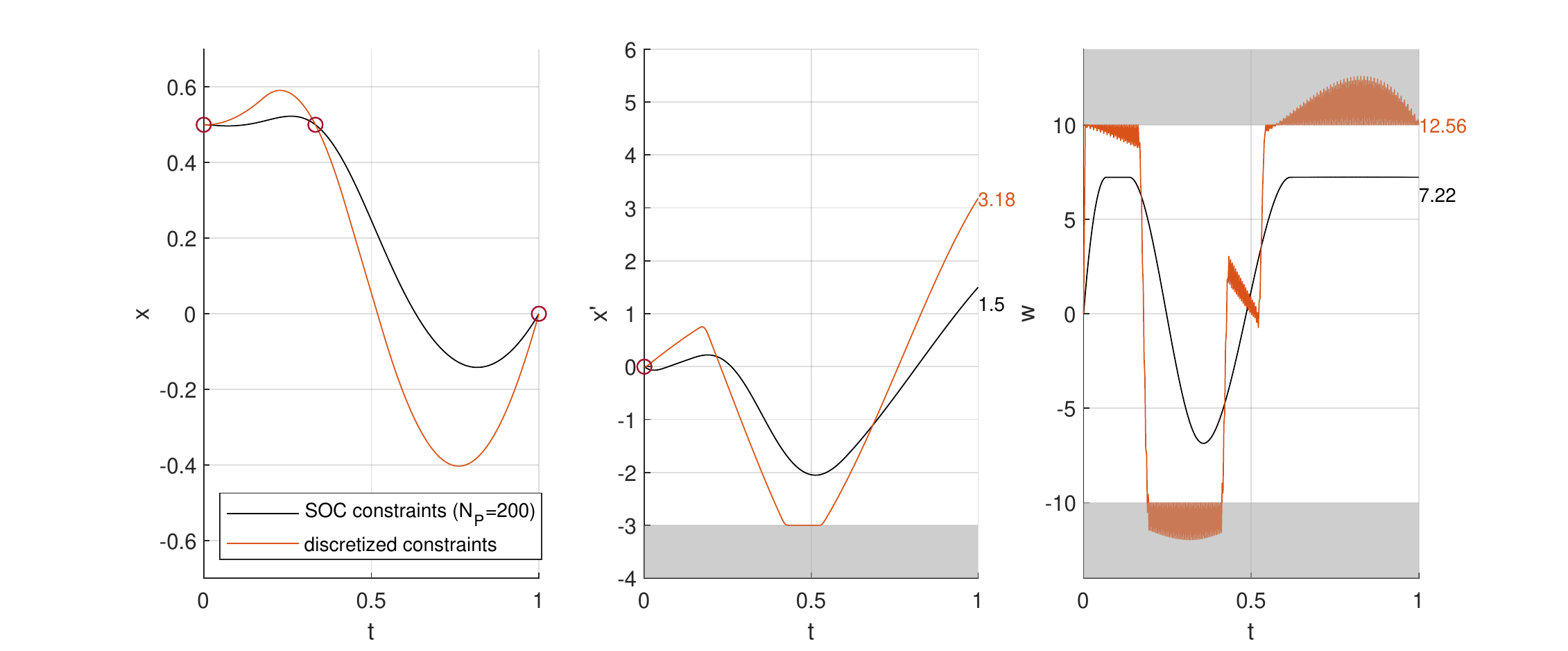}
	\includegraphics[keepaspectratio, width=\textwidth,height=\textheight]{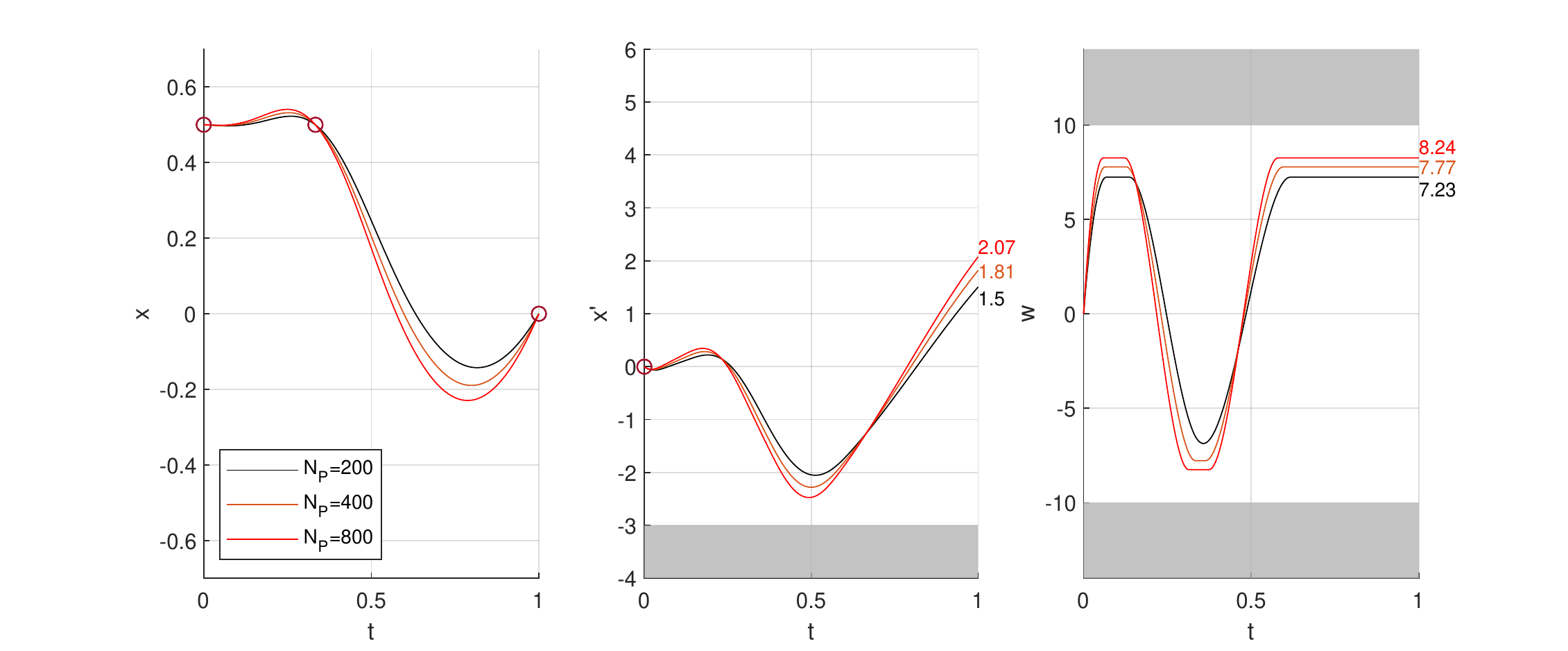}
	\includegraphics[keepaspectratio, width=\textwidth,height=\textheight]{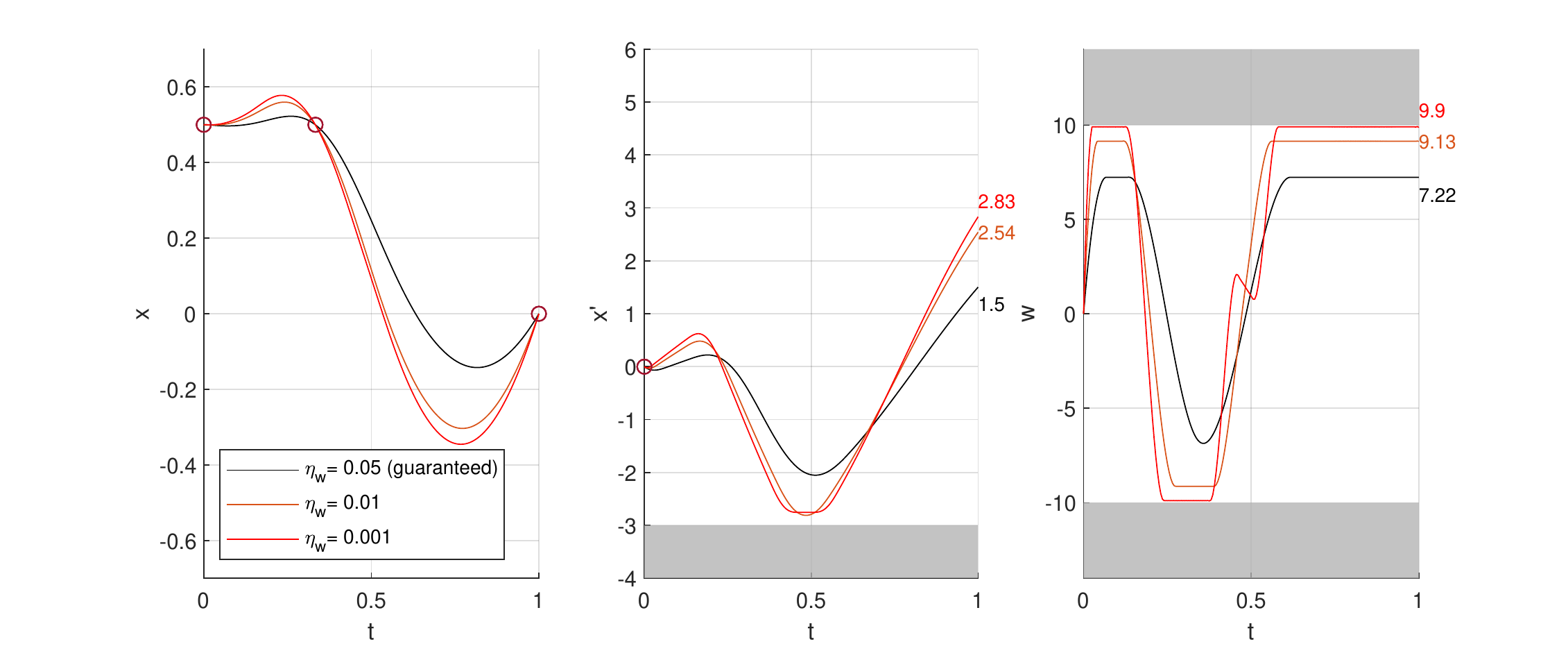}
	\caption{Optimal solutions of \eqref{opt-pendulum} for varying $N_P$ and $\eta_w$. The red circles indicate the equality-constrained points, the grayed areas the constraints over $[0,T]$. We report the values of $\dot{ x}(T)$ and of the maximum of $w(\cdot)$. \rv{ \textbf{Top:} Comparison of SOC constraints (guaranteed $\eta_w$) versus discretized constraints ($\eta_w=0$) for $N_P=200$.} \textbf{Center:} Comparison of SOC constraints for varying $N_P$ and guaranteed $\eta_w$. \textbf{Bottom:} Comparison of SOC constraints for varying $\eta_w$ and $N_P=200$.}\label{fig_pendulum}
\end{figure}
	\rv{The objective of \eqref{opt-pendulum} can be interpreted as trying to maximize the terminal velocity $\dot{ x}(T)$ with an $L^2$-cost over the control.} This example underlines how \eqref{opt-LQR_SOC_pathPlan} effectively allows for more possibilities than \eqref{opt-LQR_SOC}. We have both an initial, an intermediate and a final condition on $ x(\cdot)$. Since the full state $\b x:=[ x, \dot{ x}, w]$ is not fixed at the intermediate time $T/3$, the problem cannot be split into independent problems over two time intervals. We take $\b R(\cdot)\equiv \lambda_u \in\RR$ in order to have \eqref{opt-pendulum} written as \eqref{opt-LQR_SOC_pathPlan}. We identify the kernels $K_0$ and $K_1$ as in \eqref{def_K0} and \eqref{def_K1} and use a uniform grid of $N_P$ points for the two affine constraints over $\dot{ x}(t)$ and $w(t)$ and turn them into SOC constraints
	\begin{align*}
		\eta_{\dot{x}} \|\b x(\cdot)\|_K -\dot{ x} (t_{m})\leq -3, \quad \eta_{w} \|\b x(\cdot)\|_K + w(t_{m}) \leq 10, \quad \eta_{w} \|\b x(\cdot)\|_K - w(t_{m}) \leq 10,\, \forall \, m\in \iv{1}{N_P}
	\end{align*}
	where $\eta_{\dot{x}}$ and $\eta_{w}$ are defined as in \eqref{def_eta} for $\b C=[0\,\text{-}1\,0 \,; 0\,0\,1\,; 0\,0\,\text{-}1]$. For the experiment of Fig. \ref{fig_pendulum}, we take $T=1$, $\lambda_T=10^6$ and $\lambda_u=10^4$. \rv{All computations took less than 30 seconds.}
	
	\rv{We first compare SOC constraints with discretized constraints ($\eta_w=0$) for a moderate value of $N_P=200$. Whereas the SOC-constrained optimal trajectory is fairly conservative w.r.t.\ to the bounds over $w(t)$,  the optimal trajectory for discretized constraints suffers from chattering and does not satisfy the $w$-constraints. This was already hinted at in the speed limit example of the introduction: the control sticks to the $w$-constraints at the points $(t_m)_m$ but violates them in between. Hence the optimal value of $\dot{ x}(T)=3.18$ is attained only by repeated violation of the $w$-constraints.}
	We then present the consequences of changing the number of grid points $N_P$ and of lowering $\eta_{w}$ w.r.t.\ to the value of its definition. The parameter $\eta_{\dot{x}}$ is kept fixed to its nominal value as it has little influence on the optimal solutions. 
	We first investigate the effect of changing $N_P$ \rv{while keeping the guaranteed value of} $\eta_{w}$, defined as in \eqref{def_eta}. We see that the threshold applied to the $w$-constraints decreases, albeit slowly. Secondly, for $N_P=200$, we divide $\eta_{w}$ first by a factor 5, then by a factor 10. The $w$-constraints threshold drastically diminishes and, for $\eta_{w}=0.001$, the inflexion at $t=0.5$ appears, as a result of the constrained arc over $\dot{ x}(\cdot)$. \rv{In the limit case, we would recover the trajectory with discretized constraints ($\eta_w=0$)}.
	
	We conclude from this example that incorporating SOC terms in the constraints proves to be beneficial, even for $\etaArr$ chosen smaller than its nominal value. \rv{Nevertheless} increasing $N_P$ leads, \rv{as stated in Theorem \ref{thm_bounds_SOC}}, to convergence to the optimal trajectory \rv{$\bar{\b x}^{0}(\cdot)$} with affine state constraints through trajectories \rv{$\bar{\b x}^{\delta,\text{fin}}(\cdot)$} that are always both feasible and interior. Moreover the shapes of the SOC-optimal trajectories provide intuition on the times and properties of constrained arcs.

	\section*{Annex: Existence of interior trajectories}\label{sec_technical}
	We provide here conditions ensuring the existence of interior trajectories for \eqref{opt-LQR}. For any $\epsArr\in\R_+^P$, let $\A_\epsilon:=\{(t,\b x)\,|\,t\in[0,T],\, \epsArr+\b C(t) \b x\leq \b d(t)\}$ and $\A_{\epsilon,t}:=\{\b x\,|(t,\b x)\in \A_\epsilon\}$. Below, for $(t,\b x)\in\A_\epsilon$, $T_{\A_\epsilon}(t,\b x)$ denotes the contingent cone to the set $\A_\epsilon$ at point $(t,\b x)$ (see e.g.\ \citet{haddad1981monotone}).
	\begin{description}
		\item  [(H1)\label{hyp_gen_lin}]$\b A(\cdot)$ and $\b B(\cdot)$ are $\C^0$-continuous. $\b C(\cdot)$ and $\b d(\cdot)$ are $\C^1$-continuous and  $\b C(0) \b x_0 < \b d(0)$.
		\item  [(H2)\label{hyp_IP_lin}] There exists $M_u>0$ such that, for all $t\in [0,T]$ and $\b x \in R\BB_N$ satisfying $\b C(t) \b x \le \b d(t)$, with $R:=(1+\|\b x_0\|)e^{T\|\b A (\cdot)\|_{L^{\infty}(0,T)}+TM_u\|\b B (\cdot)\|_{L^{\infty}(0,T)}},$ there exists $\b u_{t,x}\in M_u \BB_M$ such that
		\begin{equation}\label{eq_IP_normal_lin}
		\hspace{-1cm}\forall\, i \in I_{t,x}:=\{i\,|\, \b c_i(t)^\top \b x = d_i(t) \}, \; \b C'_i(t)^\top \b x - d'_i(t) + \b c_i(t)^\top(\b A(t) \b x + \b B(t) \b u_{t,x})<0.
		\end{equation}
	\end{description} 
	The inward-pointing condition \ref{hyp_IP_lin} is a geometrical assumption on the boundary of the constraints. In particular, \ref{hyp_IP_lin} implies that the constraint set is non degenerate, i.e.\ $\A_{0,t}$ is the closure of its interior at all times $t$.
\begin{Lemma}[Existence of interior trajectories]\label{lem_IPC_cndViab}
	Under Assumptions \ref{hyp_gen_lin} and \ref{hyp_IP_lin}, the following properties are satisfied
	\begin{enumerate}[label=\roman*), ref=.\roman*]
		\item\label{item_normal_tangent} there exists $\epsArr_0>0$, $M_v>0$, $\xi>0$, and $\eta>0$ such that for all $\epsArr\le\epsArr_0$ and all $(t,\b x)\in (\partial \A_\epsilon +(0,\eta\BB_N)) \cap \A_\epsilon\cap ([0,T]\times (R-1)\BB_N)$, there exists $\b u_{t,x}\in M_u\BB_M$ such that $\b v=\b A(t) \b x + \b B(t) \b u_{t,x} \in M_v \BB_N$ and 
		\begin{equation}\label{eq_IPC_tangent_lin}
		\b y+\delta(\b v+\xi\BB_N)\subset\A_{\epsilon, t+\delta}
		\end{equation}
		for all $\delta\in[0,\xi]$ and all $\b y\in \b x+\xi\BB_N$ such that $\b y\in\A_{\epsilon, t}$. Hence $\b v\in T_{\A_\epsilon}(t,\b x)$, and
		\item\label{item_viab_strict} there exists a trajectory in $\Sx$ satisfying strictly the affine constraints, as well as the initial condition.
	\end{enumerate}	 
\end{Lemma}
\rvS{The proof below is composed of two parts. The first part translates an inward-pointing condition written in normal form \eqref{eq_IP_normal_lin} to its tangent form \eqref{eq_IPC_tangent_lin}. This allows in the second part to invoke a viability argument to derive the existence of a strictly feasible trajectory. This approach may seem remindful of a result by \citet{Soner1986}. However, unlike \citet{Soner1986}, we do not need to explicitly construct a feasible trajectory approximating a reference one, nor derive an estimate on their distance. Moreover the proof leverages the linear structure of the dynamics and of the constraints, allowing for time-varying constraints which \citet{Soner1986} did not consider.}\\

\noindent\textbf{Proof:} 
\noindent \emph{i)} Define $\b h(t,\b x):=\b C(t) \b x - \b d(t)$ and $\bm{\xi}(t,\b x,\b u):=\b C'(t) \b x - \b d'(t) + \b C(t)(\b A(t) \b x + \b B(t) \b u)$.\\

\noindent\textit{Step 1}: We claim that there exists $\eArr>0$, $\rho>0$, and $M_v>0$ such that for all $t\in [0,T]$ and $\b x \in R\BB_N$ satisfying $\b h(t,\b x) \le 0$ with some active constraints (i.e.\ $I_{t,x}\neq  \emptyset$),  there exists $\b u_{t,x}\in \R^M$ such that $\b v=\b A(t) \b x + \b B(t) \b u\in M_v \BB_N$ and
\begin{equation}\label{eq_IP_normal_cst_lin}
\forall\, i \in I^{e}_{t,x}:=\{i\,|\, h_i(t,\b x) > -e_i\}, \;\xi_i(t,\b x,\b u_{t,x})\le-\rho.
\end{equation}
Let $\b u_{t,x}$ as in \eqref{eq_IP_normal_lin}, and set 
$$\rho_{t,x}:= -\max_{i \in I_{t,x}} \xi_i(t,\b x,\b u_{t,x})>0 \quad \quad e_{t,x}:= -\max_{i \notin I_{t,x}} h_i(t,\b x,\b u_{t,x})>0.$$ 
Since $\b \xi(\cdot,\cdot,\b u_{t,x})$ and $\b h(\cdot,\cdot)$ are continuous, we can find $\Delta_{t,x}>0$ such that 
\begin{gather*}
\sup_{\substack{\delta\in[t-\Delta_{t,x},t+\Delta_{t,x}]\cap[0,T],\\ \b w\in\BB_N}} \Vert\b h(t,\b x) -\b h(t+\delta,\b x+\Delta_{t,x} \b w)\Vert_\infty \le \frac{e_{t,x}}{2}\\
\sup_{\substack{\delta\in[t-\Delta_{t,x},t+\Delta_{t,x}]\cap[0,T],\\ \b w\in\BB_N,\, i\in I_{t,x}}} \vert \xi_i(t,\b x,\b u_{t,x}) - \xi_i(t+\delta,\b x+\Delta_{t,x} \b w,\b u_{t,x})\vert \le \frac{\rho_{t,x}}{2}.
\end{gather*}
This implies that the index set of active constraints does not increase in size for points in the open set $\Omega_{t,x}:=(t,\b x)+\Delta_{t,x} ]-1,1[\times\mathring{\BB}_N$, denoting by $\mathring{\BB}_N$ the open unit ball of $\RR^N$. Moreover, for any point in $\Omega_{t,x}$, $\b u_{t,x}$ satisfies our claim for the constants $\rho_{t,x}/2$ and $(e_{t,x}/2)\,\b 1_N$. Since we are considering a compact set, we can select a finite number of $(t_j,\b x_j)_{j\in\iv{1}{J}}$ such that 
\begin{align*}
([0,T]\times R\BB_N)\cap \A_0  \subset \bigcup_{j\in\iv{1}{J}}\Omega_{t_j,\b x_j}.
\end{align*}
To conclude, simply take\footnote{If $\A_0$ is bounded, take $M_u:= 1+ \max_{j \in \iv{1}{J}} \|\b u_{t_j,\b x_j}\|$.}
\begin{gather*}
\rho:= \min_{j \in \iv{1}{J}} \frac{\rho_{t_j,\b x_j}}{2},\quad \quad \b e:= \left(\min_{j \in \iv{1}{J}} \frac{e_{t_j,\b x_j}}{2}\right) \,\b 1_N, \quad \quad
M_v:= R \|\b A(\cdot)\|_{L^{\infty}(0,T)} + M_u\|\b B(\cdot)\|_{L^{\infty}(0,T)}.
\end{gather*}
\noindent\textit{Step 2}: Thanks to the uniform constants of the normal form \eqref{eq_IP_normal_cst_lin}, we derive the constants of the tangent form \eqref{eq_IPC_tangent_lin}. Let us choose $\xi,\eta\in]0,1]$ such that for all $(t,\b x)\in \partial \A_0\cap ([0,T]\times R\BB_N)$, setting $\b v_{t,x}=\b A(t) \b x + \b B(t) \b u_{t,x} \in M_v \BB_N$, for
$$\b e^{t,x,\delta}_{ \bm{\alpha},\bm{\beta},\bm{\gamma}}:=(\b C(t+\delta) -\b C(t))(\b x + \eta \bm{\gamma} + \xi \bm{\alpha}) + \delta \b C(t+\delta) (\b v_{t,x}+\eta \b A(t)\bm{\gamma}+ \xi \bm{\beta})+\b d(t) -\b d(t+\delta),$$
we have
\begin{align}
\frac{\b e}{2}> \sup_{\substack{\delta\in[t-\xi,t+\xi]\cap[0,T]\\ \bm{\alpha},\bm{\beta},\bm{\gamma}\in\BB_N}} \b e^{t,x,\delta}_{ \bm{\alpha},\bm{\beta},\bm{\gamma}} \quad, \quad \frac{-\delta\rho}{2} > \sup_{\substack{\delta\in[t-\xi,t+\xi]\cap[0,T]\\ \bm{\alpha},\bm{\beta},\bm{\gamma}\in\BB_N,\, i\in I^e_{t,x}}}  (\b e^{t,x,\delta}_{ \bm{\alpha},\bm{\beta},\bm{\gamma}})_i. \label{eq_tang_nori,mPC}
\end{align}
The first inequality can be derived from the $\C^1$-smoothness of $\b C(\cdot)$ and $\b d(\cdot)$, and from the fact that both $\b x$, $\b v_{t,x}$, $\b A(\cdot)$ and $\b C(\cdot)$ are uniformly bounded. The second inequality stems from \eqref{eq_IP_normal_cst_lin}, as $\lim_{\delta\rightarrow 0^+}\frac{\b e^{t,x,\delta}_{0,0,0}}{\delta}=\b \xi(t,\b x,\b u_{t,x})$.
The two inequalities of \eqref{eq_tang_nori,mPC} state that, for $\b y:=\b x + \eta \bm{\gamma} + \xi \bm{\alpha}$, if $\b y\in \A_{0,t}$ then $$\b y+\delta(\b v+\xi\BB_N)\subset\A_{0, t+\delta} \text{ with } \b v:=\b v_{t,x}+\eta \b A(t)\bm{\gamma}=\b A(t) (\b x+\eta \bm{\gamma}) + \b B(t) \b u_{t,x}.$$
\noindent\textit{Step 3}: There just remains to extend the result of Step 2 to the perturbed constraint sets $\A_\epsilon$. For any $t\in[0,T]$, let $\b w_t$ be an eigenvector of the largest eigenvalue $\mu_t$ of $\b C(t) \b C(t)^T$ satisfying $\|\b w_t\|_\infty=\max_{i\in\iv{1}{P}} w_{t,i}=1$. Let $\bar{\mu}:=\min_{t\in[0,T]} \mu_t>0$. Take any $\epsArr_0>0$ such that
\begin{align*}
\sup_{t\in[0,T]} \frac{2 \|\epsArr_0\|_\infty \|\b C(t)^\top \b w_t\|}{\bar{\mu}}\le \frac{\eta}{2}.
\end{align*}
Let $\epsArr\le\epsArr_0$, take $(t,\b x)\in \partial \A_\epsilon \cap ([0,T]\times (R-1)\BB_N)$. Then 
$$ \exists\, i\in\iv{1}{P},\; \b c_i(t)^\top\left(\b x + \frac{2 \|\epsArr_0\|_\infty \b C(t)^\top \b w_t}{\bar{\mu}}\right) = d_i(t) -\epsArr_i+   \frac{2\mu_t}{\bar{\mu}} \|\epsArr_0\|_\infty >0.$$
Since $\eta\le 1$, one can therefore find $\breve{\b x}\in R\BB_N \cap \partial \A_{0,t}$ such that $\|\b x - \breve{\b x}\|=d_{\A_{0,t}}(\b x)\leq \eta/2$, to which the conclusions of Step 2 apply. So the set considered in \eqref{eq_IPC_tangent_lin} at $\b x$, taking as constants $\xi/2$ and $\eta/2$, is a subset of the one defined by $\breve{\b x}$, $\xi$ and $\eta$. Since only differences of $\b d(\cdot)$ appeared in Step 2 and that $\b x\in \partial \A_{\epsilon,t}$, adding $\epsArr$ has no effect on the computations, so it follows that \eqref{eq_IPC_tangent_lin} is satisfied. For $\b y= \b x$, by definition of $T_{\A_\epsilon}(t,\b x)$, we have that $(1,\b v)\in T_{\A_\epsilon}(t,\b x)$.\\

\noindent \emph{ii)} Let $\b F(t,\b x):=\{(1, \b A(t) \b x + \b B(t) \b u)\,|\, \b u\in M_u\BB_M\}$. Since $\b A(\cdot)$ and $\b B(\cdot)$ are $\C^0$-continuous, the bounded set-valued map $\b F(\cdot,\cdot)$ is upper semicontinuous. Let $\Omega:=\A_\epsilon\cap ([0,T]\times (R-1)\mathring{\BB}_N)$. The open ball $\mathring{\BB}_N$ being open in $\RR^N$, $\Omega$ is locally compact in $[0,T]\times \R^N$. In \emph{i)} we have shown that the local viability condition is satisfied for $\A_\epsilon$. The intersection with an open set does not add boundary points, so, for $(t,\b x)\in \Omega$, $\b F(t,\b x)\cap T_{\Omega}(t,\b x) \neq \emptyset$. We may therefore apply \citep[Theorem 1]{haddad1981monotone} which provides a trajectory $\b x^\epsilon(\cdot)$ satisfying $\b x^\epsilon(0)=\b x_0$, $\b x^{\epsilon,'}(t) \in \b F(t,\b x^{\epsilon}(t))$, and $\b h(t,\b x^{\epsilon}(t))+\epsArr \le 0$. Let $[0,t_1[$ be the maximal interval of existence of $\b x^{\epsilon}(\cdot)$. Since $\b x^{\epsilon,'}(\cdot)$ is measurable, by \citep[Theorem 2.3.13]{vinter2010optimal}, we can find some measurable $\b u^\epsilon(\cdot)$ with values bounded by $M_u$, s.t.\ $\b x^{\epsilon,'}(\cdot)=\b A(t) \b x + \b B(t) \b u^\epsilon(t)$ a.e.\ in $[0,T]$ . The dynamics being sublinear as $\b A(\cdot)$,  $\b B(\cdot)$, and  $\b u^\epsilon(\cdot)$ are bounded, $\b x^\epsilon(\cdot)$ can be continuously extended to $(t_1,\b x^\epsilon(t_1))\in \Omega$. So $t_1\geq T$, otherwise the viability condition would allow to extend $\b x^\epsilon(\cdot)$ beyond $t_1$. Since $\b u^\epsilon(\cdot)$ is measurable and bounded, $\b u^\epsilon(\cdot)\in L^2(0,T)$. Hence $\b x^\epsilon(\cdot)\in\Sx$ satisfies the required properties.

\begin{flushright}
	$\blacksquare$
\end{flushright}

\textbf{Acknowledgments}

\rv{The author is grateful to the four anonymous referees for their comments and their constructive suggestions that helped to improve this article.} The author also expresses his recognition to Nicolas Petit for the enlightening discussions when writing down these investigations.

	\bibliographystyle{apalike}
	\bibliography{LQ_HardCons} 

\begin{thebibliography}{}

\bibitem[Ackermann, 1985]{Ackermann1985}
Ackermann, J. (1985).
\newblock {\em Sampled-Data Control Systems}.
\newblock Springer Berlin Heidelberg.

\bibitem[Aronszajn, 1950]{aronszajn50theory}
Aronszajn, N. (1950).
\newblock Theory of reproducing kernels.
\newblock {\em Transactions of the American Mathematical Society}, 68:337--404.

\bibitem[Aubin-Frankowski, 2020]{aubin2020Riccati}
Aubin-Frankowski, P.-C. (2020).
\newblock {Interpreting the dual Riccati equation through the LQ reproducing
  kernel}.
\newblock {\em Comptes Rendus. Math\'ematique}.
\newblock (\url{https://arxiv.org/abs/2012.12940}).

\bibitem[Aubin-Frankowski et~al., 2020]{aubin2020ifac}
Aubin-Frankowski, P.-C., Petit, N., and Szab{\'o}, Z. (2020).
\newblock Kernel regression for vehicle trajectory reconstruction under speed
  and inter-vehicular distance constraints.
\newblock In {\em IFAC World Congress}, volume (to appear).

\bibitem[Aubin-Frankowski and Szab{\'o}, 2020a]{aubin2020hard_SDP}
Aubin-Frankowski, P.-C. and Szab{\'o}, Z. (2020a).
\newblock Handling hard affine {SDP} shape constraints in {RKHSs}.
\newblock Technical report.
\newblock (\url{https://arxiv.org/abs/2101.01519}).

\bibitem[Aubin-Frankowski and Szab{\'o}, 2020b]{aubin2020hard_nips}
Aubin-Frankowski, P.-C. and Szab{\'o}, Z. (2020b).
\newblock Hard shape-constrained kernel machines.
\newblock In {\em Advances in Neural Information Processing Systems (NIPS)}.

\bibitem[Bertalan et~al., 2019]{bertalan2019learning}
Bertalan, T., Dietrich, F., Mezić, I., and Kevrekidis, I.~G. (2019).
\newblock On learning {Hamiltonian} systems from data.
\newblock {\em Chaos: An Interdisciplinary Journal of Nonlinear Science},
  29(12):121107--1 -- 121107--9.

\bibitem[Bonnans and Ram{\'{\i}}rez, 2005]{Bonnans2005}
Bonnans, J.~F. and Ram{\'{\i}}rez, H. (2005).
\newblock Perturbation analysis of second-order cone programming problems.
\newblock {\em Mathematical Programming}, 104(2-3):205--227.

\bibitem[Bourdin and Trélat, 2017]{bourdin2017linearquadratic}
Bourdin, L. and Trélat, E. (2017).
\newblock Linear–quadratic optimal sampled-data control problems:
  {Convergence} result and {Riccati} theory.
\newblock {\em Automatica}, 79:273--281.

\bibitem[Burachik et~al., 2014]{burachik2014duality}
Burachik, R.~S., Kaya, C.~Y., and Majeed, S.~N. (2014).
\newblock A {Duality} {Approach} for {Solving} {Control}-{Constrained}
  {Linear}-{Quadratic} {Optimal} {Control} {Problems}.
\newblock {\em SIAM Journal on Control and Optimization}, 52(3):1423--1456.

\bibitem[Chaplais et~al., 2011]{chaplais2011design}
Chaplais, F., Malisani, P., and Petit, N. (2011).
\newblock Design of penalty functions for optimal control of linear dynamical
  systems under state and input constraints.
\newblock In {\em Conference on Decision and Control (CDC)}, pages 6697--6704.

\bibitem[Chiuso and Pillonetto, 2019]{chiuso2019system}
Chiuso, A. and Pillonetto, G. (2019).
\newblock System {Identification}: {A} {Machine} {Learning} {Perspective}.
\newblock {\em Annual Review of Control, Robotics, and Autonomous Systems},
  2(1):281--304.

\bibitem[Dower et~al., 2019]{dower2019barriers}
Dower, P.~M., McEneaney, W.~M., and Cantoni, M. (2019).
\newblock Game representations for state constrained continuous time linear
  regulator problems.
\newblock {\em arXiv:1904.05552}.
\newblock http://arxiv.org/abs/1904.05552.

\bibitem[Fujii and Kawahara, 2019]{Fujii2019}
Fujii, K. and Kawahara, Y. (2019).
\newblock {Dynamic mode decomposition in vector-valued reproducing kernel
  Hilbert spaces for extracting dynamical structure among observables}.
\newblock {\em Neural Networks}, 117:94--103.

\bibitem[Fujioka and Kano, 2013]{fujioka2013control}
Fujioka, H. and Kano, H. (2013).
\newblock Control theoretic {B}-spline smoothing with constraints on
  derivatives.
\newblock pages 2115--2120.

\bibitem[Gerdts and Hüpping, 2012]{gerdts2012virtual}
Gerdts, M. and Hüpping, B. (2012).
\newblock Virtual control regularization of state constrained linear quadratic
  optimal control problems.
\newblock {\em Computational Optimization and Applications}, 51(2):867--882.

\bibitem[Giannakis et~al., 2019]{giannakis2019reproducing}
Giannakis, D., Das, S., and Slawinska, J. (2019).
\newblock {Reproducing kernel Hilbert space compactification of unitary
  evolution groups}.
\newblock http://arxiv.org/abs/1808.01515.

\bibitem[Grüne and Guglielmi, 2018]{grune2018MPC}
Grüne, L. and Guglielmi, R. (2018).
\newblock Turnpike properties and strict dissipativity for discrete time linear
  quadratic optimal control problems.
\newblock {\em SIAM Journal on Control and Optimization}, 56(2):1282--1302.

\bibitem[Haddad, 1981]{haddad1981monotone}
Haddad, G. (1981).
\newblock Monotone trajectories of differential inclusions and functional
  differential inclusions with memory.
\newblock {\em Israel Journal of Mathematics}, 39(1-2):83--100.

\bibitem[Hartl et~al., 1995]{hartl_survey_1995}
Hartl, R.~F., Sethi, S.~P., and Vickson, R.~G. (1995).
\newblock A {Survey} of the {Maximum} {Principles} for {Optimal} {Control}
  {Problems} with {State} {Constraints}.
\newblock {\em SIAM Review}, 37(2):181--218.

\bibitem[Heckman, 2012]{heckman2012theory}
Heckman, N. (2012).
\newblock The theory and application of penalized methods or {Reproducing}
  {Kernel} {Hilbert} {Spaces} made easy.
\newblock {\em Statistics Surveys}, 6(0):113--141.

\bibitem[Hermant, 2009]{hermant2009stab}
Hermant, A. (2009).
\newblock Stability analysis of optimal control problems with a second-order
  state constraint.
\newblock {\em SIAM Journal on Optimization}, 20(1):104--129.

\bibitem[{Kailath}, 1971]{kailath1971RKHS}
{Kailath}, T. (1971).
\newblock {RKHS approach to detection and estimation problems--I: Deterministic
  signals in Gaussian noise}.
\newblock {\em IEEE Transactions on Information Theory}, 17(5):530--549.

\bibitem[{Kano} and {Fujioka}, 2018]{kano2018Bsplines}
{Kano}, H. and {Fujioka}, H. (2018).
\newblock B-spline trajectory planning with curvature constraint.
\newblock In {\em 2018 Annual American Control Conference (ACC)}, pages
  1963--1968.

\bibitem[Kojima and Morari, 2004]{kojima2004lq}
Kojima, A. and Morari, M. (2004).
\newblock {LQ} control for constrained continuous-time systems.
\newblock {\em Automatica}, 40(7):1143--1155.

\bibitem[Magnus~Egerstedt, 2009]{egerstedt2009control}
Magnus~Egerstedt, C.~M. (2009).
\newblock {\em Control Theoretic Splines: Optimal Control, Statistics, and Path
  Planning (Princeton Series in Applied Mathematics)}.
\newblock Princeton University Press.

\bibitem[Marco et~al., 2017]{marco2017design}
Marco, A., Hennig, P., Schaal, S., and Trimpe, S. (2017).
\newblock On the design of {LQR} kernels for efficient controller learning.
\newblock In {\em {Conference} on {Decision} and {Control} ({CDC})}, pages
  5193--5200. IEEE.

\bibitem[Mattingley and Boyd, 2012]{mattingley2012CVX}
Mattingley, J. and Boyd, S. (2012).
\newblock {CVXGEN}: a code generator for embedded convex optimization.
\newblock {\em Optimization and Engineering}, 13(1):1--27.

\bibitem[{Mercy} et~al., 2016]{mercy2016MotionPlan}
{Mercy}, T., {Van Loock}, W., and {Pipeleers}, G. (2016).
\newblock Real-time motion planning in the presence of moving obstacles.
\newblock In {\em European Control Conference (ECC)}, pages 1586--1591.

\bibitem[Micheli and Glaunès, 2014]{micheli_matrix-valued_2014}
Micheli, M. and Glaunès, J.~A. (2014).
\newblock Matrix-valued kernels for shape deformation analysis.
\newblock {\em Geometry, Imaging and Computing}, 1(1):57--139.

\bibitem[Parzen, 1970]{parzen1970statis}
Parzen, E. (1970).
\newblock Statistical inference on time series by {RKHS} methods.
\newblock In {\em Proceedings 12th Biennial Seminar, Canadian Mathematical
  Congress}.

\bibitem[Pillonetto et~al., 2014]{pillonetto2014kernel}
Pillonetto, G., Dinuzzo, F., Chen, T., Nicolao, G.~D., and Ljung, L. (2014).
\newblock Kernel methods in system identification, machine learning and
  function estimation: {A} survey.
\newblock {\em Automatica}, 50(3):657 -- 682.

\bibitem[Rosenfeld et~al., 2019]{rosenfeld2019dynamic}
Rosenfeld, J.~A., Kamalapurkar, R., Gruss, L.~F., and Johnson, T.~T. (2019).
\newblock {Dynamic Mode Decomposition for Continuous Time Systems with the
  Liouville Operator}.
\newblock http://arxiv.org/abs/1910.03977.

\bibitem[Schölkopf et~al., 2001]{Scholkopf2001}
Schölkopf, B., Herbrich, R., and Smola, A.~J. (2001).
\newblock A generalized representer theorem.
\newblock In {\em Computational Learning Theory (CoLT)}, pages 416--426.
  Springer Berlin Heidelberg.

\bibitem[Sch{\"o}lkopf and Smola, 2002]{scholkopf02learning}
Sch{\"o}lkopf, B. and Smola, A. (2002).
\newblock {\em Learning with Kernels: Support Vector Machines, Regularization,
  Optimization, and Beyond}.
\newblock MIT Press.

\bibitem[Schwartz, 1964]{Schwartz1964}
Schwartz, L. (1964).
\newblock Sous-espaces hilbertiens d’espaces vectoriels topologiques et
  noyaux associés (noyaux reproduisants).
\newblock {\em Journal d’Analyse Mathématique}, 13:115--256.

\bibitem[Singh et~al., 2018]{singh2018kernelStab}
Singh, S., Sindhwani, V., Slotine, J.-J., and Pavone, M. (2018).
\newblock Learning stabilizable dynamical systems via control contraction
  metrics.
\newblock In {\em Workshop on Algorithmic Foundations of Robotics (WAFR)}.
\newblock https://arxiv.org/abs/1808.00113.

\bibitem[Soner, 1986]{Soner1986}
Soner, H.~M. (1986).
\newblock Optimal control with state-space constraint i.
\newblock {\em {SIAM} Journal on Control and Optimization}, 24(3):552--561.

\bibitem[Sootla et~al., 2018]{Sootla2018}
Sootla, A., Mauroy, A., and Ernst, D. (2018).
\newblock {Optimal control formulation of pulse-based control using Koopman
  operator}.
\newblock {\em Automatica}, 91:217--224.

\bibitem[Speyer and Jacobson, 2010]{speyer2010primer}
Speyer, J.~L. and Jacobson, D.~H. (2010).
\newblock {\em Primer on Optimal Control Theory}.
\newblock Society for Industrial and Applied Mathematics, USA.

\bibitem[Steinke and Schölkopf, 2008]{steinke2008kernels}
Steinke, F. and Schölkopf, B. (2008).
\newblock Kernels, regularization and differential equations.
\newblock {\em Pattern Recognition}, 41(11):3271--3286.

\bibitem[van {Keulen}, 2020]{Keulen2020MPC}
van {Keulen}, T. (2020).
\newblock Solution for the continuous-time infinite-horizon linear quadratic
  regulator subject to scalar state constraints.
\newblock {\em IEEE Control Systems Letters}, 4(1):133--138.

\bibitem[Van~Loan, 1978]{VanLoan1978computingExpM}
Van~Loan, C. (1978).
\newblock Computing integrals involving the matrix exponential.
\newblock {\em IEEE Transactions on Automatic Control}, 23(3):395--404.

\bibitem[Vinter, 1990]{vinter2010optimal}
Vinter, R. (1990).
\newblock {\em Optimal Control}.
\newblock Birkhäuser, Basel.

\bibitem[Wahba, 1990]{wahba90spline}
Wahba, G. (1990).
\newblock {\em Spline Models for Observational Data}.
\newblock SIAM, CBMS-NSF Regional Conference Series in Applied Mathematics.

\bibitem[Williams et~al., 2015]{OWilliams2015}
Williams, M.~O., , Rowley, C.~W., Kevrekidis, I.~G., and and (2015).
\newblock {A kernel-based method for data-driven Koopman spectral analysis}.
\newblock {\em Journal of Computational Dynamics}, 2(2):247--265.

\end{thebibliography}
\end{document}